\documentclass[10pt]{article}

\usepackage{graphicx,hyperref}
\usepackage{amsmath}
\usepackage[version=4]{mhchem}
\usepackage{siunitx}
\usepackage{textcomp}
\usepackage{longtable,tabularx}
\usepackage{bm}
\usepackage{cite}
\usepackage{float}
\usepackage{pstricks}
\usepackage{amsfonts}
\usepackage{latexsym}
\usepackage{enumerate}
\usepackage{epsfig}
\usepackage{subfig}
\usepackage{framed,fancybox}
\usepackage{bbm}
\usepackage{multirow}
\usepackage{dashrule}
\usepackage[T1]{fontenc}
\usepackage{threeparttable}
\usepackage{rotating}
\usepackage{amssymb,amscd}
\usepackage[ruled,vlined]{algorithm2e}
\usepackage[bottom]{footmisc}
\allowdisplaybreaks
\makeatletter
\def\myVCENTER#1{\vcenter{\hbox{$\m@th#1$}}}
\makeatother

\long\def\symbolfootnote[#1]#2{\begingroup\def\thefootnote{\fnsymbol{footnote}}\footnote[#1]{#2}\endgroup}

\definecolor{shadecolor}{gray}{0.99}
{\endMakeFramed}

\newenvironment{shadedframe}{%
 \MakeFramed {\FrameRestore}}
{\endMakeFramed}

\long\def\symbolfootnote[#1]#2{\begingroup\def\thefootnote{\fnsymbol{footnote}}\footnote[#1]{#2}\endgroup}
\def\qed{\hfill{$\vcenter{\hrule height1pt \hbox{\vrule width1pt height5pt
    \kern5pt \vrule width1pt} \hrule height1pt}$} \medskip}
\newcommand{\m}[1]{{\bf{#1}}}

\newcommand{\tr}{^{\sf T}}
\newcommand{\C}[1]{{\cal {#1}}}

\newcolumntype{C}{ >{\centering\arraybackslash} m{4cm} }
\newcolumntype{D}{ >{\centering\arraybackslash} m{1cm} }
\newcommand{\R}{\mathbb R}

\newcommand{\J}{\mathcal{J}}
\newcommand{\M}{\mathcal{M}}

\newcommand{\dt}{\text{ d}}
\newcommand{\crb}[1]{\{ #1 \}}

\newcommand{\deriv}[2]{\dfrac{\text{d} #1}{\text{d} #2}}

\title{\textbf{Method for Solving State-Path Constrained \\ Optimal Control Problems Using \\ Adaptive Radau Collocation}}

\author{Cale~A.~Byczkowski\footnote{Ph.D.~Candidate, Department of Mechanical and Aerospace Engineering.  Email: cale.byczkowski@ufl.edu.} \\ Anil V.~Rao\footnote{Professor, Department of Mechanical and Aerospace Engineering, University Term Professor.  Associate Fellow, AIAA.  E-mail:  anilvrao@ufl.edu.  Corresponding Author} \\ \date{\today} \\ {\em University of Florida} \\ {\em Gainesville, FL 32611--6250}}

\date{}
\begin{document}

\renewcommand{\baselinestretch}{1}
\normalsize\normalfont 
\maketitle

\begin{abstract}
A new method is developed for accurately approximating the solution to state-variable inequality path constrained optimal control problems using a multiple-domain adaptive Legendre-Gauss-Radau collocation method. The method consists of the following parts. First, a structure detection method is developed to estimate switch times in the activation and deactivation of state-variable inequality path constraints. Second, using the detected structure, the domain is partitioned into multiple-domains where each domain corresponds to either a constrained or an unconstrained segment. Furthermore, additional decision variables are introduced in the multiple-domain formulation, where these additional decision variables represent the switch times of the detected active state-variable inequality path constraints. Within a constrained domain, the path constraint is differentiated with respect to the independent variable until the control appears explicitly, and this derivative is set to zero along the constrained arc while all preceding derivatives are set to zero at the start of the constrained arc. The time derivatives of the active state-variable inequality path constraints are computed using automatic differentiation and the properties of the chain rule. The method is demonstrated on two problems, the first being a benchmark optimal control problem which has a known analytical solution and the second being a challenging problem from the field of aerospace engineering in which there is no known analytical solution. When compared against previously developed adaptive Legendre-Gauss-Radau methods, the results show that the method developed in this paper is capable of computing accurate solutions to problems whose solution contain active state-variable inequality path constraints.
\end{abstract}

\renewcommand{\baselinestretch}{1.5}
\normalsize\normalfont

\pagebreak
\section{Introduction}
An optimal control problem is one where it is desired to optimize a specified performance index while satisfying dynamic constraints, path constraints, and boundary conditions.  Since most optimal control problems do not have analytic solutions, numerical methods must be developed and used to approximate the solution of an optimal control problem.  Numerical methods for optimal control fall into two broad categories: indirect and direct methods. In an indirect method, the first-order optimality conditions are derived using the calculus of variations and the original optimal control problem is converted to a Hamilton boundary-value problem (HBPV). The HBVP is then solved using either single-shooting, multiple-shooting, or indirect collocation. In a direct method, either the control or the state and control are approximated at a set of support points, and the original optimal control problem is transcribed to a nonlinear programming problem (NLP). The NLP can then solved using well-developed software such as \textit{SNOPT}\cite{Gill2002}, \textit{IPOPT}\cite{Biegler2008}, or \textit{KNITRO}\cite{Byrd2007}.

A particular class of direct methods (known as \textit{direct collocation}) has been extensively used for computing the numerical approximation to a continuous time optimal control problem. In a direct collocation method (also referred to as an implicit simulation method), both the state and control are parameterized, and the constraints of the continuous time optimal control problem are enforced at a set of specially chosen \textit{collocation} points. In the past decade, a great deal of research has gone into the development of \textit{Gaussian quadrature orthogonal collocation methods}. In a Gaussian quadrature orthogonal collocation method, the state is approximated using a basis of orthogonal Lagrange polynomials where the specially chosen set of collocation points is associated with a Gaussian quadrature. The most well developed Gaussian quadrature methods employ either Legendre-Gauss (LG) points~\cite{Benson2006,Rao2010}, Legendre-Gauss-Radau (LGR) points~\cite{Kameswaran1,Garg2011June,Garg2011April,Garg2010}, or Legendre-Gauss-Lobatto (LGL) points~\cite{Elnagar1995}. It has been shown that when the solution for the optimal control takes on a smooth form, Gauss quadrature methods are capable of converging to highly accurate solutions at an exponential rate~\cite{Hou2012}. Conversely, when the solution is nonsmooth or singular, which is often the case when constraints are enforced on the state and/or control, advanced collocation methods tend to require a large amount of computation time and often return poor or no solution for the problem of interest. As a result, standard direct collocation methods often return suboptimal or even unfeasible solutions. 

In the past decade, a class of so called \textit{hp}-adaptive mesh refinement methods have been developed that can handle discontinuities in the solution and compute highly accurate solutions \cite{Darby:2011a,Darby:2011b,Patterson:2015,Liu:2017,Liu2018}. Though, it is often the case that an \textit{hp}-adaptive method will place an unnecessary amount of collocation points and mesh intervals around a discontinuity, resulting in a larger amount of computation. A similar behavior is exhibited when a state-variable inequality constraint is present in the problem, thus motivating the development of a method that can properly partition a problem such that the rapid convergence of a Gaussian quadrature collocation method can still be exploited without suffering when a state-variable inequality constraint is present.

Increased-complexity optimal control problems are often subject to path constraints.  An important class of path-constrained optimal control problems that poses challenges computationally are problems with {\em state-variable inequality constraints} (SVICs), henceforth referred to as a state-path constrained optimal control problem (SPOC). That is, a SPOC is an optimal control problem with inequality constraints that are only a function of the state and possibly time. The added complexity of solving a SPOC typically stems from not knowing when and where the inequality constraints are active on the optimal solution.  Given the constraint is only a function of the state (and possibly time), no information on the optimal control is readily available when the constraint is active.  Since the early 1960s, a number of researchers have delved into establishing additional necessary conditions for optimality when such constraints are present.  This research has led to a variety of different sets of optimality conditions, mainly stemming from three different approaches often referred to as the method of \textit{direct adjoining} \cite{Jacobson1971,Chang1962,Norris1973,McIntyre1967}, \textit{indirect adjoining} \cite{Bryson1963,Pontryagin,Dreyfus1962,Leitmann1971,Kreindler1982}, and \textit{indirect adjoining with continuous multipliers} \cite{Russak1,Russak2,Russak3,Gollan1980,Hartl1995}. 

In the direct adjoining approach, the SVICs (in their original form) are adjoined to the Hamiltonion of the system through the use of Lagrange multipliers. Next, the calculus of variations is used to derive first-order optimality conditions which result in a set of "jump conditions" for the optimal costate at the beginning and end of any constrained arcs, leading to possible discontinuities in the costate. As for the indirect adjoining approach, obtaining the necessary conditions for optimality involves adjoining to the Hamiltonian the lowest time derivative of the SVIC which is explicit in the control. Using this approach, a set of necessary "tangency conditions" are enforced at either the entrance, and/or exit of the contrained arc, which also lead to discontinuities in the costate \cite{Bryson1963}. Lastly, the approach of indirect adjoining with continuous multipliers involves defining a new costate variabale so as to "subtract" the discontinuity, leading to a continuous costate \cite{Hartl1995}.  A small class of problems exist for which an analytical solution can be derived when applying the optimality conditions obtained from any of the aforementioned approaches, and for this reason, research for developing numerical methods to approximate the solution to a SPOC has been a focus for decades.

As alluded to earlier, from a numerical point of view, additional challenges exist when trying to approximate the solution to inequality-constrained optimal control problems. For example, when the inequality constraints are purely a function of the control, the optimal control profile may contain discontinuities and singular arcs, which can be difficult to approximate numerically. In the case that inequality constraints are purely a function of the state variables, the problem leads to a system consisting of \textit{high-index} differential-algebraic equation (DAE) constraints when the constraint is active, which are computationally challenging to solve with traditional numerical integration methods \cite{biehn2000,Feehery1998,Betts2002}. Moreover, activations of the SVICs may induce nonsmooth behavior in the state and control solutions while also introducing discontinuities in the costates. This research is motivated by the importance of having a numerical method capable of computing an accurate approximation to problems whose solutions may be nonsmooth due to the existence of SVICs.

Several numerical methods have been developed for solving a SPOC, all similarly falling within the category of either a direct, indirect, or a combination of both; all of which utilize the various forms of the additional necessary conditions for inequality constraints. The most common methods that appear within the indirect approach are so called penalty function techniques \cite{Liu2016,Mall2020,Heidrich2021} and transformation techniques \cite{Jacobson1969,Graichen2009,Fabien2013}. In a penalty function method, the original constrained problem is approximated by augmenting the inequality constraint to the cost function through the use of a penalty term. The approximated problem, which is now unconstrained, is iteratively solved until it approaches the solution of the original problem. In a transformation method, the dimension of the state is increased through the inclusion of parametric variables (also known as ``slack variables"). If the order of the inequality constraint is $q$, then the dimension of the state increases by $q$, and the additional state variables are defined such that the inequality constraint is satisfied at every point throughout the trajectory. Both penalty term and transcription methods benefit from the fact that no a priori knowledge of when the constraint becomes active is necessary; however, both methods suffer computationally as each need to employ some form of iterative process to obtain an accurate approximate of the solution. 

As for handling SVICs in a direct approach, the state variable constraint (or an appropriate time derivative) are directly adjoined to the problem\cite{Dreyfus1962,Denham1964,Feehery1998,Kameswaran2008}. Specifically, the state path constraint and its $q-1$ time derivatives are enforced as equality constraints at the entry of the constrained arc, while the $q^{th}$ time derivative is enforced as a mixed state/control equality constraint. A pitfall in using a direct approach stems from the fact that the number and location of constraint activations must be guessed and/or known a priori. 

The contribution of this work is a new method for solving state-path constrained optimal control problems.  The method described in this paper consists of the following parts.  First, a structure detection algorithm is developed to detect estimated activation and deactivation times of the SVICs. Second, based on the results of the detection algorithm, the multiple-domain reformulation of the LGR collocation method developed in Ref.~\cite{Pager2022} is used to partition the problem into subdomains where the detected constraints are either deemed to be active or inactive. Similar to Ref.~\cite{Pager2022}, the multiple-domain formulation partitions the time horizon such that the estimated activation and deactivation times are included as additional decision variables in the optimization process. Third, the method automatically enforces the additional necessary conditions following from the approach of indirect adjoining. Specifically, a method is developed to algorithmically obtain the higher time derivatives of the active SVICs in order to form the necessary tangency conditions which are enforced at the beginning of constrained domains while the highest such time derivative of the active SVICs that is explicit in a control variable is enforced as an equality constraint during the constrained domains. Lastly, the method developed in this paper utilizes adaptive LGR quadrature to suggest refined meshes as well as allow for the ability to re-perform structure detection. This adds an additional level of robustness to the method to allow for the capability to accurately detect the constraint structure.

%The contributions of this work are as follows. First, a structure detection algorithm is developed to obtain estimates of activation and deactivation times of any active SVICs. Second, based on the detection results, the method automatically partitions the original time domain into multiple subdomains where the constraints are either set to be active or inactive. Third, the method does not require any a priori knowledge of the structure of the active SVICs. Fourth, the method includes the detected activation and deactivation times within the optimization decision vector, allowing for to the exact locations to be optimized. Fifth, the method includes additional necessary conditions for optimality within the direct method formulation. Specifically, the approach of breaking up the original time domain into multiple subdomains allows for the use of additional necessary conditions to be extended to direct methods.

The remainder of the paper is organized as follows. Section~\ref{section:bolza} introduces an optimal control problem in Bolza form. Section~\ref{section:MDLGR} describes the multiple-domain formulation of Legendre-Gauss-Radau collocation used within this work to formulate the nonlinear programming problem. Section~\ref{section:svic} provides an overview of the additional necessary conditions used within this work for handling SVICs. Section~\ref{section:method} details the method developed in this work for solving state-path constrained optimal control problems. Section~\ref{section:examples} demonstrates the method on two examples. Section~\ref{section:limitations} discusses limitations of the method developed in this work. Finally, Section~\ref{section:conclusions} provides conclusions on this work.

%%%%%%%%%%%%%%%%%%%%%%%%%%%%%%% SECTION 2 %%%%%%%%%%%%%%%%%%%%%%%%%%%%%%%%%%%%%%%%%%%%%%%%%%%%%%%
\section{Bolza Optimal Control Problem}\label{section:bolza}
Consider the following general optimal control problem in Bolza form. Determine the state, $\m{y}(t) \in \R^{n_y}$, the control, $\m{u}(t) \in \R^{n_u}$, the initial time, $t_0$, and the terminal time, $t_f$, on the time interval, $t \in [t_0, t_f]$, (where $n_y$ and $n_u$ represent the number of states and controls, respectively) that minimizes the cost functional
\begin{equation}\label{eq:optprb1_cost}
\J = \M(\m{y}(t_0),t_0,\m{y}(t_f),t_f) + \int_{t_0}^{t_f} \mathcal{L}(\m{y}(t),\m{u}(t),t) \dt t,
\end{equation}
subject to the dynamic constraints
\begin{equation}\label{eq:optprb1_dyn}
\dfrac{\dt\m{y}}{\dt t} = \m{f}(\m{y}(t),\m{u}(t),t),
\end{equation}
the inequality path constraints 
\begin{equation}\label{eq:optprb1_ineq}
\m{c}_{\text{min}} \leq \m{c}(\m{y}(t),\m{u}(t),t) \leq \m{c}_{\text{max}},
\end{equation}
and the boundary conditions
\begin{equation}\label{eq:optprb1_bc}
\m{b}_{\text{min}} \leq \m{b}(\m{y}(t_0),t_0,\m{y}(t_f),t_f) \leq \m{b}_{\text{max}}.
\end{equation}
The functions $\M$, $\mathcal{L}$, $\m{f}$, $\m{b}$, $\m{c}$ are defined by the following mappings
\begin{equation*}
\begin{array}{lclcl}
\M &:& \R^{n_y} \times \R \times \R^{n_y} \times \R &\rightarrow& \R,   \\
\mathcal{L} &:& \R^{n_y} \times \R^{n_u} \times \R &\rightarrow& \R,   \\
\m{f} &:& \R^{n_y} \times \R^{n_u} \times \R &\rightarrow& \R^{n_y}, \\
\m{b} &:& \R^{n_y} \times \R \times \R^{n_y} \times \R &\rightarrow& \R^{n_b}, \\
\m{c} &:& \R^{n_y} \times \R^{n_u} \times \R &\rightarrow& \R^{n_c},
\end{array}
\end{equation*}
where $n_b$ and $n_c$ are the number of specified boundary conditions and inequality path constraints, respectively. From a numerical standpoint, and without loss of generality, it is useful to map the original time interval, $t \in [t_0,t_f]$, of the Bolza optimal control problem onto the time interval $\tau \in [-1,+1]$. This is done with the use of the following affine transformation
\begin{equation}
t(\tau,t_0,t_f) = \dfrac{t_f - t_0}{2}\tau + \dfrac{t_f - t_0}{2}.
\end{equation}
The Bolza optimal control problem of Eqs.(\ref{eq:optprb1_cost})-(\ref{eq:optprb1_bc}) can now be defined in terms of $\tau$ as follows. Determine the state, $\m{y}(\tau) \in \R^{n_y}$, the control, $\m{u}(\tau) \in \R^{n_u}$, the initial time, $t_0$, and the terminal time, $t_f$, on the time interval, $\tau \in [-1, +1]$, that minimizes the cost functional
 \begin{equation}\label{eq:optprb2_cost}
\J = \M(\m{y}(-1),t_0,\m{y}(+1),t_f) + \dfrac{t_f - t_0}{2}\int_{-1}^{+1} \mathcal{L}(\m{y}(\tau),\m{u}(\tau),t(\tau,t_0,t_f)) \dt \tau,
\end{equation}
subject to the dynamic constraints
\begin{equation}\label{eq:optprb2_dyn}
\dfrac{\dt\m{y}}{\dt \tau} = \dfrac{t_f - t_0}{2} \m{f}(\m{y}(\tau),\m{u}(\tau),t(\tau,t_0,t_f)),
\end{equation}
the inequality path constraints 
\begin{equation}\label{eq:optprb2_ineq}
\m{c}_{\text{min}} \leq \m{c}(\m{y}(\tau),\m{u}(\tau),t(\tau,t_0,t_f)) \leq \m{c}_{\text{max}},
\end{equation}
and the boundary conditions
\begin{equation}\label{eq:optprb2_bc}
\m{b}_{\text{min}} \leq \m{b}(\m{y}(-1),t_0,\m{y}(+1),t_f) \leq \m{b}_{\text{max}}.
\end{equation}

%%%%%%%%%%%%%%%%%%%%%%%%%%%%%%% SECTION 3 %%%%%%%%%%%%%%%%%%%%%%%%%%%%%%%%%%%%%%%%%%%%%%%%%%%%%%%
\section{Multiple-Domain Legendre-Gauss-Radau Collocation}\label{section:MDLGR}
While any valid collocation scheme can be used, the method developed in this paper employs Legendre-Gauss-Radau (LGR) collocation to discretize the continuous Bolza optimal control problem. The LGR collocation method is chosen due to its ability to converge to a solution at an exponential rate when the solution is smooth \cite{Hou2012}. However, for the problems studied in this research, the solution may be nonsmooth due to active state-constrained arcs. Consequently, in order to be able to correctly enforce additional necessary conditions and handle potential nonsmooth behavior, the method developed in this research utilizes a multiple-domain LGR collocation method (similar to Pager, et. al.~\cite{Pager2022}) described below. For a visual representation of the domain partitioning, see Figure~\ref{fig:structureDecomp}.
\subsection{Multiple-Domain Formulation}
Let the horizon $t \in [t_0,t_f]$ of the Bolza optimal control problem be partitioned into $D$ domains, $\C{P}_d = [t_s^{\crb{d-1}},t_s^{\crb{d}}]  \subseteq [t_0,t_f]$, $d = 1,...,D$, such that
\begin{equation}
\begin{array}{lcl}
\bigcup_{d=1}^{D} \C{P}_d &=& [t_0,t_f], \\[5pt]
\bigcap_{d=1}^{D} \C{P}_d &=& \{t_s^{\crb{1}},...,t_s^{\crb{D-1}}\},
\end{array}
\end{equation}
where $t_s^{\crb{d}}$, are the \textit{domain interface variables}, $t_s^{\crb{0}}=t_0$, $t_s^{\crb{D}}=t_f$, are the initial and final times, respectively, $d$ is the domain index, and $s$ signifies a switch between a constrained or unconstrained arc. Specifically, for the method developed in this research, $t_s^{\crb{d}}$ for $d \neq \{0,D\}$ correspond to switches in activation and deactivation of the state-variable inequality constraint, whereas in Pager, et. al. \cite{Pager2022}, these variables serve as switch points in bang-bang and/or singular control structures. It is noted that the domain interface variables become additional decision variables in the resulting nonlinear programming problem (NLP) and are \textit{not} collocation points.

Similar to the formulation presented before, each domain $\C{P}_d = [t_s^{\crb{d-1}},t_s^{\crb{d}}]$ is mapped to $\tau \in [-1,+1]$ using the affine transformation
\begin{equation}\label{eq:affine_map}
t(\tau,t_s^{\crb{d-1}},t_s^{\crb{d}}) = \dfrac{t_s^{\crb{d}} - t_s^{\crb{d-1}}}{2}\tau + \dfrac{t_s^{\crb{d}} - t_s^{\crb{d-1}}}{2}  
\end{equation}
Suppose now that the time interval $\tau \in [-1,+1]$ on each domain is divided into $K$ mesh intervals, $\C{S}_k = [T_{k-1},T_{k}] \subseteq [-1,+1]$, for $k \in \crb{1,...,K}$, such that
\begin{equation}
\begin{array}{lcl}
\bigcup_{k=1}^{K} \C{S}_k &=& [-1,+1], \\[5pt]
\bigcap_{k=1}^{K} \C{S}_k &=& \{T_1,...,T_{K-1}\},
\end{array}
\end{equation}
and $-1 = T_0 < T_1 < \cdots < T_{K-1} < T_{K} = +1$.
\subsection{Legendre-Gauss-Radau Collocation}
Once the multiple-domain formulation is complete, LGR collocation is used to discretize the continuous time problem on each mesh interval. Specifically, for each mesh interval, the LGR points are defined on $[T_{k-1},T_{k}] \subseteq [-1,+1]$, $k \in \crb{1,...,K}$. The state of the continuous time optimal control problem is then approximated within each mesh interval $\C{S}_k$ by
\begin{equation}\label{eq:state_approx}
\m{y}^{(k)}(\tau) \approx \m{Y}^{(k)}(\tau) = \sum_{j=1}^{N_k + 1}\m{Y}_j^{(k)}\ell_j^{(k)}(\tau), \quad \ell_j^{(k)}(\tau) = \prod_{\substack{i = 1 \\ i \neq j}}^{N_k+1} \dfrac{\tau - \tau_i^{(k)}}{\tau_j^{(k)}-\tau_i^{(k)}}, \quad k \in \crb{1,...,K},
\end{equation}
where $\ell_j^{(k)}(\tau)$, $j \in \crb{1,...,N_k + 1}$ is a basis of Lagrange polynomials on $S_k$, and $(\tau_1^{(k)},...,\tau_{N_k}^{(k)})$ are the $N_k$ set of LGR collocation points in the interval $[T_{k-1},T_k)$, with $\tau_{N_k+1} = T_k$ being a non-collocated support point; the first collocation point of the next interval is equal to the non-collocated point in the previous interval $\tau_1^{(k)} = T_{k-1}$. Differentiating Eq.~\eqref{eq:state_approx} with respect to $\tau$ leads to
\begin{equation}
\deriv{\m{Y}^{(k)}(\tau)}{\tau} = \sum_{j=1}^{N_k+1} \m{Y}_j^{(k)} \deriv{\ell_j^{(k)}(\tau)}{\tau}.
\end{equation}
The dynamic constraints given in Eq.~\eqref{eq:optprb2_dyn} are then approximated at the $N_k$-LGR points in mesh interval $\C{S}_k$, $k \in \crb{1,...,K}$ on domain $d \in \crb{1,...,D}$ by
\begin{equation}\label{approxdynconstr}
\sum_{j=1}^{N_k+1}D_{ij}^{(k)}\m{Y}_j^{(k)} - \dfrac{t_f - t_0}{2} \m{f}\left( \m{Y}_i^{(k)}, \m{U}_i^{(k)}, t(\tau_i^{(k)},t_s^{\crb{d-1}},t_s^{\crb{d}}) \right) = 0, \hspace{0.5em} \quad i \in \crb{1,...,N_k}
\end{equation}
where
\begin{equation}
D_{ij}^{(k)} = \dfrac{\text{d} \ell_j^{(k)}(\tau)}{\text{d} \tau}, \hspace{0.5em} \quad i \in \crb{1,...,N_k}, \quad j \in \crb{1,...,N_k+1},
\end{equation}
are the elements of the $N_k \times (N_k + 1)$ Legendre-Gauss-Radau \textit{differentiation matrix} in the mesh interval $\C{S}_k$, and $\m{U}_i^{(k)}$ is the control parameterization at the $i$-th collocation point in mesh interval $\C{S}_k$. Continuity in the state across mesh intervals $\mathcal{S}_{k-1},\mathcal{S}_{k}$ and time domains $\C{P}_{d-1},\C{P}_d$ is achieved by using the same variable to represent $\m{Y}_{N_{k-1}+1}^{k-1} = \m{Y}^{(k)}_1$ and $\m{Y}_{N^{\crb{d-1}}+1}^{\crb{d-1}} = \m{Y}^{\crb{d}}_1$, respectively. Lastly, it is noted that $N^{\crb{d}}$ is the total number of collocation points in the time domain $\C{P}_d$ and can be computed by,
\begin{equation}
N^{\crb{d}} = \sum_{k=1}^{K^{\crb{d}}} N_k^{\crb{d}}
\end{equation}
where $K^{\crb{d}}$ is the total number of mesh intervals in time domain $\C{P}_d$, $d \in \crb{1,...,D}$. 

\subsection{Nonlinear Programming Problem}\label{section:NLP}
The multiple-domain Legendre-Gauss-Radau (LGR) discretization of the continuous time multiple-domain Bolza optimal control problem results in a large sparse nonlinear programming (NLP) which is given as follows. Minimize the objective functional
\begin{equation}\label{eq:NLP_cost}
\C{J} = \C{M} \left( \m{Y}_1^{\crb{1}}, t_0, \m{Y}_{N^{\crb{D}}+1}^{\crb{D}}, t_f \right) + \sum_{d=1}^{D}\dfrac{t_s^{\crb{d}}-t_s^{\crb{d-1}}}{2} \left[\m{w}^{\crb{d}}\right]\tr \m{L}^{\crb{d}},
\end{equation} 
subject to the dynamic constraints
\begin{equation}\label{eq:NLP_dynconstr}
\m{\Delta}^{\crb{d}} = \m{D}^{\crb{d}}\m{Y}^{\crb{d}} - \dfrac{t_s^{\crb{d}}-t_s^{\crb{d-1}}}{2} \m{F}^{\crb{d}} = \m{0}\hspace{0.5em}, \quad d \in \crb{1,...,D},
\end{equation}
the inequality path constraints
\begin{equation}\label{eq:NLP_ineq}
\m{c}_{\text{min}} \leq \m{C}^{\crb{d}} \leq \m{c}_{\text{max}}, \hspace{0.5em} \quad d \in \crb{1,...,D},
\end{equation}
the boundary conditions
\begin{equation}
\m{b}_{\text{min}} \leq \m{b}\left( \m{Y}_1^{\crb{1}}, t_0, \m{Y}_{N^{\crb{D}}+1}^{\crb{D}}, t_f \right) \leq  \m{b}_{\text{max}},
\end{equation}
and the continuity conditions
\begin{equation}\label{eq:contconds}
\m{Y}_{N^{\crb{d-1}}+1}^{\crb{d-1}} = \m{Y}^{\crb{d}}_1,
\end{equation}
where $\m{D}^{\crb{d}} \in \R^{N^{\crb{d}} \times (N^{\crb{d}}+1)}$ is the LGR differentiation matrix in the time domain $\C{P}_d$, $d \in \crb{1,...,D}$, and $\m{w}^{\crb{d}} \in \R^{N^{\crb{d}} \times 1}$ are the LGR weights associated with the approximation of the integral of the cost function via gaussian quadrature. The decision variables of the NLP given by Eqs.~\eqref{eq:NLP_cost}-\eqref{eq:contconds} are the elements of the state matrix $\m{Y}^{\crb{d}} \in \R^{(N^{\crb{d}}+1) \times n_y}$, the elements of the control matrix $\m{U}^{\crb{d}} \in \R^{N^{\crb{d}} \times n_u}$, the initial time $t_s^{\crb{0}}=t_0$, the terminal time $t_s^{\crb{D}}=t_f$, and the domain interface variables $t_s^{\crb{d}}$ for $d \in \crb{1,...,D-1}$. Note, the continuity conditions in Eq.~\eqref{eq:contconds} are implicitly enforced by using the same decision variable in the NLP for $\m{Y}_{N^{\crb{d-1}}+1}^{\crb{d-1}}$ and $\m{Y}^{\crb{d}}_1$. The state matrix, $\m{Y}^{\crb{d}} \in \R^{(N^{\crb{d}}+1) \times n_y}$, and the control matrix, $\m{U}^{\crb{d}} \in \R^{N^{\crb{d}} \times n_u}$, in time domain $\C{P}_d$, $d \in \crb{1,...,D}$ are formed as
\begin{equation}\label{eq:state_control_approx}
\m{Y}^{\crb{d}} = \begin{bmatrix}
\m{Y}^{\crb{d}}_1 \\ \vdots \\ \m{Y}^{\crb{d}}_{N^{\crb{d}}+1}
\end{bmatrix} \hspace{0.5em} \text{, and } 
\m{U}^{\crb{d}} = \begin{bmatrix}
\m{U}^{\crb{d}}_1 \\ \vdots \\ \m{U}^{\crb{d}}_{N^{\crb{d}}}
\end{bmatrix},
\end{equation}
respectively.  
Lastly, the elements $\m{F}^{\crb{d}}, \m{C}^{\crb{d}}, \m{L}^{\crb{d}}$, in Eqs.~\eqref{eq:NLP_cost}-\eqref{eq:NLP_ineq} for $d \in \crb{1,...,D}$ are defined as follows
\begin{equation}\label{eq:aprrox_eqs}
\begin{array}{lclcl}
\m{F}^{\crb{d}} &=& \begin{bmatrix}
\m{f}\left( \m{Y}_1^{\crb{d}}, \m{U}_1^{\crb{d}}, t_1^{\crb{d}} \right) \\
\vdots \\
\m{f}\left( \m{Y}_{N^{\crb{d}}}^{\crb{d}}, \m{U}_{N^{\crb{d}}}^{\crb{d}}, t_{N^{\crb{d}}}^{\crb{d}} \right)
\end{bmatrix} &\in& \R^{N^{\crb{d}}\times n_y}, \vspace{1em} \\ 
\end{array}
\end{equation}
\begin{equation}
\begin{array}{lclcl}
\m{C}^{\crb{d}} &=& \begin{bmatrix}
\m{c}\left( \m{Y}_1^{\crb{d}}, \m{U}_1^{\crb{d}}, t_1^{\crb{d}} \right) \\
\vdots \\
\m{c}\left( \m{Y}_{N^{\crb{d}}}^{\crb{d}}, \m{U}_{N^{\crb{d}}}^{\crb{d}}, t_{N^{\crb{d}}}^{\crb{d}} \right)
\end{bmatrix} &\in& \R^{N^{\crb{d}}\times n_c}, \vspace{1em} \\
\end{array}
\end{equation}
\begin{equation}
\begin{array}{lclcl}
\m{L}^{\crb{d}} &=& \begin{bmatrix}
\mathcal{L} \left( \m{Y}_1^{\crb{d}}, \m{U}_1^{\crb{d}}, t_1^{\crb{d}} \right) \\
\vdots \\
\mathcal{L} \left( \m{Y}_{N^{\crb{d}}}^{\crb{d}}, \m{U}_{N^{\crb{d}}}^{\crb{d}}, t_{N^{\crb{d}}}^{\crb{d}} \right)
\end{bmatrix} &\in& \R^{N^{\crb{d}}\times 1},  \\
\end{array}
\end{equation}
where $n_y$ is the number of state components and $n_u$ is the number of control components in the problem.
\subsection{Approximation of Solution Error}\label{section:solution_error}
Suppose the nonlinear program of Eqs.~\eqref{eq:NLP_cost}-\eqref{eq:contconds} has been solved on a mesh, $\mathcal{S}_k$, $k \in \crb{1,...,K}$ on domain $\mathcal{P}_d$, $d \in \crb{1,...,D}$ with $N_k$ collocation points in mesh interval $\mathcal{S}_k$. An estimate of the discretization error on the current mesh must be obtained in order to assess the accuracy of the solution. The method developed in this paper directly uses the maximum estimated relative error in the state solution as a condition for terminating the algorithm. The process used for estimating the relative error is identical to Ref. \cite{Patterson:2015}. The process is summarized as follows. \par
The objective is to approximate the error in the state at a set of $M_k = N_k + 1$ LGR points $( \hat{\tau}_1^{(k)}, ..., \hat{\tau}_{M_k}^{(k)} )$, where $\hat{\tau}_1^{(k)} = \tau_1^{(k)} = T_{k-1}$ and $\hat{\tau}_{M_k + 1}^{(k)} = T_k$. Let the state approximation at the $M_k$ LGR points be denoted by $( \m{Y}(\hat{\tau}_1^{(k)}),..., \m{Y}(\hat{\tau}_{M_k}^{(k)}))$. Next, let the control be approximated in mesh interval $\mathcal{S}_k$ with a Lagrange polynomial
\begin{equation}
\m{U}^{(k)}(\tau) = \sum_{j=1}^{N_k} \m{U}_j^{(k)}\hat{\ell}_j^{(k)}(\tau), \quad \hat{\ell}_j^{(k)}(\tau) = \prod_{\substack{l = 1 \\ l \neq j}}^{N_k} \dfrac{\tau - \tau_l^{(k)}}{\tau_j^{(k)} - \tau_l^{(k)}},
\end{equation}
and let the control approximation at $\hat{\tau}_i^{(k)}$ be denoted $\m{U}(\hat{\tau}_i^{(k)})$ for $i = 1,...,M_k$, then the value of the right-hand side of the dynamics at $(\m{Y}(\hat{\tau}_i^{(k)}),\m{U}(\hat{\tau}_i^{(k)}),\hat{\tau}_i^{(k)})$ is used to construct an improved approximation of the state. Let $\hat{\m{Y}}^{(k)}$ be a polynomial of degree at most $M_k$ that is defined on the mesh interval $\C{S}_k$. If the derivative of $\hat{\m{Y}}^{(k)}$ matches the dynamics at each LGR quadrature points $\hat{\tau}_i^{(k)}$, $i = 1,...,M_k$, then we have
\begin{equation}
\begin{array}{lcl}
\hat{\m{Y}}^{(k)}(\hat{\tau}_{j+1}^{(k)}) = \m{Y}^{(k)}(\hat{\tau}_1^{(k)}) + \dfrac{t_f - t_0}{2} \sum\limits_{l = 1}^{M_k} \hat{I}_{jl}^{(k)} \m{f}(\m{Y}^{(k)}(\hat{\tau}_l^{(k)}),\m{U}^{(k)}(\hat{\tau}_l^{(k)}),t(\hat{\tau}_l^{(k)},t_0,t_f)), \left\lbrace \begin{array}{c}
j = 1,...,M_k \\ k = 1,...,K
\end{array} \right\rbrace
\end{array}
\end{equation}
where $\hat{I}_{jl}^{(k)}$, $j,l = 1,...,M_k$, is the $M_k \times M_k$ LGR integration matrix corresponding to the LGR points defined by $( \hat{\tau}_1^{(k)}, ..., \hat{\tau}_{M_k}^{(k)} )$.
Comparing the interpolated values, $\m{Y}(\hat{\tau}_l^{(k)})$, $l = 1,...,M_k + 1$, with the integrated values $\hat{\m{Y}}(\hat{\tau}_l^{(k)})$, $l = 1,...,M_k + 1$, the \textit{relative error} in the $i^{th}$ component of the state at $( \hat{\tau}_1^{(k)}, ..., \hat{\tau}_{M_k + 1}^{(k)} )$ is defined as
\begin{equation}
\begin{array}{lcl}
e_i^{(k)}(\hat{\tau}_l^{(k)}) &=& \dfrac{\left| \hat{Y}^{(k)}_i(\hat{\tau}_l^{(k)}) - Y^{(k)}_i(\hat{\tau}_l^{(k)}) \right|}{\displaystyle1 + \max_{\stackrel{j \in\{1,\ldots,N_k+1\}}{k \in\{1,\ldots,K\}}} \left| Y^{(k)}_i(\tau_j^{(k)}) \right|}, \quad \left\lbrace \begin{array}{c}
l = 1,...,M_k + 1 \\ k = 1,...,K \\i = 1,...,n_y
\end{array} \right\rbrace.
\end{array}
\end{equation}
The \textit{maximum relative error} on domain $\C{P}_d$ in mesh interval $\C{S}_k$ is then defined as
\begin{equation}\label{eq:max_mesh_err}
\begin{array}{lcl}
e_{\max}^{(k)} &=& \displaystyle \max_{\stackrel{i \in\{1,\ldots,n_y\}}{l \in\{1,\ldots,M_k+1\}}} e_i^{(k)}(\hat{\tau}_l^{(k)})
\end{array}
\end{equation}

%%%%%%%%%%%%%%%%%%%%%%%%%%%%%%% SECTION 4 %%%%%%%%%%%%%%%%%%%%%%%%%%%%%%%%%%%%%%%%%%%%%%%%%%%%%%%
\section{State-Variable Inequality Constrained Optimal Control}\label{section:svic}
Suppose now that the inequality path constraint given by Eq.~\eqref{eq:optprb1_ineq} takes on the following form
\begin{equation}\label{eq:pathconstr3}
\m{c}(\m{y}(t),t) \leq \m{0}.
\end{equation}
In order to determine a condition for the optimal control when the inequality in Eq.~\eqref{eq:pathconstr3} holds, successive time derivatives of Eq.~\eqref{eq:pathconstr3} are taken, and the dynamic constraints $\m{f}(\m{y}(t),\m{u}(t),t)$ are substituted for $\text{d}\m{y}(t)/\text{d}t$ until an expression explicit in $\m{u}(t)$ is obtained. This approach is commonly referred to as the method of indirect adjoining \cite{Bryson1963,Pontryagin,Dreyfus1962,Leitmann1971,Kreindler1982}. If $q$ successive time derivatives are required, then the state-variable inequality path constraints (SVICs) given by Eq.~\eqref{eq:pathconstr3} are said to be of $q$th-order. The Hamiltonian for the system becomes 
\begin{equation}\label{eq:adjoined_Hamiltonian}
\C{H} = \C{L} + \bm{\lambda} \tr \m{f} + \bm{\mu}\tr \m{c}^{(q)},
\end{equation}
where $\bm{\lambda}$ are the costate variables, $\bm{\mu}$ are the Lagrange multipliers associated with the SVICs that satisfy the following necessary conditions
\begin{align}
\bm{\mu} &\geq \m{0}, \text{ on the constraint boundary }(\m{c}^{(q)} = \m{0}), \\
\bm{\mu} &= \m{0}, \text{ off the constraint boundary }(\m{c} < \m{0}),
\end{align}
and
\begin{equation}
\m{c}^{(q)} = \dfrac{\text{d}^{(q)}\m{c}}{\text{d}t^{(q)}}.
\end{equation}
The first-order optimality conditions then become
\begin{equation}
\dot{\bm{\lambda}}^{\sf T}  = -\dfrac{\partial \C{H}}{\partial \m{y}} = 
\begin{cases} 
-\dfrac{\partial \mathcal{L}}{\partial \m{y}} - \bm{\lambda} ^{\sf T}  \dfrac{\partial \m{f}}{\partial \m{y}} - \bm{\mu} \dfrac{\partial \m{c}^{(q)}}{\partial \m{y}}, & \m{c}=\m{0} \vspace{1em}\\
-\dfrac{\partial \mathcal{L}}{\partial \m{y}} - \bm{\lambda} ^{\sf T} \dfrac{\partial \m{f}}{\partial \m{y}}, & \m{c}<\m{0}
\end{cases}
\end{equation}
In order to compute a finite control to keep the system on the constraint boundary, the path entering the constraint must satisfy the following \textit{tangency} constraints \cite{Bryson1963}
\begin{equation}\label{eq:tangencyconstraints}
\m{N}(\m{y}(t_a),t_a) \equiv \begin{bmatrix}
\m{c}(\m{y},t_a) \\[5pt] 
\m{c}^{(1)}(\m{y},t_a) \\[5pt] 
\vdots \\[5pt] 
\m{c}^{(q-1)}(\m{y},t_a)
\end{bmatrix}
= \m{0},
\end{equation}
where $t_a$ is the time at which the SVICs become active. It has been shown in \cite{biehn2000,Feehery1998,Betts2002} that when a constraint is active, the index of the original system of differential algebraic equations (DAEs) increases. In a direct method, enforcing Eq.~\eqref{eq:tangencyconstraints} at the start \textit{or} end of a constrained arc is analogous to performing index reduction.  Specifically, the tangency constraints provide a set of necessary conditions that reduce the index of the constrained system of DAEs to the index of the unconstrained system. Although direct methods generally perform well when path constraints are present in the problem, it is not possible to enforce the tangency conditions at the correct instances throughout the trajectory without prior knowledge of when the SVICs become active or inactive. The challenge of being able to apply the tangency constraints within a direct method motivates the development of the method presented in this work, which is described in the next section.

%%%%%%%%%%%%%%%%%%%%%%%%%%%%%%% SECTION 5 %%%%%%%%%%%%%%%%%%%%%%%%%%%%%%%%%%%%%%%%%%%%%%%%%%%%%%%
\section{Method for State-Path Constrained Optimal Control Problems}\label{section:method}

In this section the method for solving state-path constrained optimal control problems (SPOC) is presented, and will be called the SPOC method. Recall, an optimal control problem with inequality path constraints only a function of the state and possibly time will be referred henceforth as a SPOC. The SPOC method is comprised of two stages. The first stage, described in further detail in Section \ref{section:structure_detect_decomp}, employs a structure detection algorithm to identify activation and deactivation (A/D) times in each state-variable inequality path constraint (SVIC). Next, the detected A/D times are treated as domain interface variables and are used to partition the domain of the independent variable into subdomains such that each subdomain consists of segments where the detected constraint is either active or inactive. The domain interface variables are then treated as additional decision variables in the nonlinear programming problem (NLP) in order to optimize the location(s) for when the detected constraint is active. The next stage of the method described in Section \ref{section:constraints_and_refinement} enforces additional necessary conditions for optimality such that the partitioned problem is properly constrained. Specifically, time derivatives of the detected constraints are used to reduced the order of the differential algebraic system of equations and properly constrain the allowable control search space. The final step in the method, as shown in Section \ref{section:procedure}, is to determine if it is necessary to perform mesh refinement and/or structure detection and decomposition to improve the accuracy of the solution.

\subsection{Structure Detection and Decomposition}\label{section:structure_detect_decomp}
The discussion to follow will assume a single SVIC is present in the problem formulation, but the method developed in this research can be generalized to handle multiple SVICs. Assume now that the optimal control problem formulated in Sect. \ref{section:bolza} subject to a SVIC of the form defined in Sect. \ref{section:svic} has been transcribed into a nonlinear programming problem (NLP) using the multiple-domain LGR collocation method provided in Sect. \ref{section:MDLGR} with $D = 1$ (that is, a single domain is used). The solution of the NLP results in a discrete approximation of the state and control, as given by Eq. (\ref{eq:state_control_approx}). Assume further that the maximum relative error does not meet the prescribed mesh refinement accuracy tolerance. As a result, mesh refinement is required which simultaneously enables the decomposition of the independent variable into subdomains where the SVIC is either active or inactive. The subdomains are formulated using structure detection as described below. \par
Structure detection obtains estimates of activation and deactivation (A/D) times in the SVIC on a provided mesh (in this case, the initial mesh) in order to classify subdomains as either constrained or unconstrained. In this work, only the approximation of the SVIC (given by Eq.~\eqref{eq:aprrox_eqs}) is used in the detection process. Structure detection begins by applying the method of Sect. \ref{section:structure_detect} to identify and estimate any A/D times in the SVIC. After A/D locations have been estimated, the method of Sect. \ref{sec:structure_decomp} takes the estimated A/D times and determines the classification of each domain as either constrained or unconstrained. The structure detection process developed in this section is always applied on the initial mesh. Structure detection is only applied recursively if the maximum constraint violations on the current mesh are larger than a user specified constraint satisfaction tolerance ($\epsilon_c$). It is important to note that structure detection is the crux of the method developed in this work and the ability to recursively perform structure detection will ensure an accurate estimate of the true A/D times can be obtained.

\subsubsection{Structure Detection}\label{section:structure_detect}
The overall approach used in this research aims to automatically decompose the original inequality constrained optimal control problem into a multiple-domain optimal control problem consisting of state constrained and unconstrained domains. Partitioning the problem into multiple domains provides the ability to apply the higher order necessary tangency conditions at the beginning of constrained domains. The multiple-domain formulation also allows the resulting mixed state-control equality constraint (lowest time derivative of the path constraint that first contains the control) to be enforced throughout constrained domains. Applying the higher order necessary conditions reduces the high-index system of differential algebraic equations (DAEs), and provides the ability to accurately approximate the control while any SVICs are active.

The first step in partitioning the problem is to obtain an estimate of the domain interface variables, which are analogous to activation/deactivation (A/D) times of any SVICs. Before presenting the structure detection method for estimating the domain interface variables, the following assumption is made in the development of the algorithm:
\begin{itemize}
\item[] \textbf{Assumption 1:} The SVICs are not active at the start of the trajectory: $\m{c}(\m{y}(t_0),t_0) < \m{0}$.
\end{itemize}
The above assumption is in place to be able to inspect neighboring collocation points to determine the activation of a constraint arc and properly enforce the necessary tangency conditions at the start of a constrained domain. It is noted that in Jacobson, et. al.~\cite{Jacobson1971}, a set of new necessary conditions were derived which show that, in certain cases, for SVICs of \textit{odd} order ($q > 1$), the state may \textit{not} be constrained over a non-zero interval of time, but instead only touch the constraint boundary at a single point. This occurrence is often referred to as a \textit{touch point}. It is mentioned here that the structure detection method developed in this research is intended for SVICs that contain constrained arcs over non-zero intervals of time. That is, the potential existence of touch points are not studied in this research.  \par 
The SPOC method begins by computing an approximate solution to the optimal control problem on a single domain $\C{P}_1 = [t_0,t_f]$ using standard Legendre-Gauss-Radau (LGR) collocation on a coarse fixed mesh consisting of $K$ mesh intervals $S_k,$ with $N_k$ number of collocation points within each mesh interval. Next, let $\m{Y}(\tau_j^{(k)})$, $j \in \crb{1,...,N_k+1}$, $k \in \crb{1,...,K}$, be the approximated state solution obtained on the initial domain at the $j^{th}$ collocation point in the $k^{th}$ mesh interval, then the following relative difference is computed at every collocation point plus the final non-collocated point of the final mesh interval for every SVIC in the problem
\begin{equation}\label{eq:svic_rel_diff}
\delta c_i (\tau_j^{(k)}) =  \dfrac{\left| c_i(\m{Y}(\tau_j^{(k)}),\tau_j^{(k)}) - c_{i,\max/\min} \right|}{1 + |c_{i,\max/\min}|}, \hspace{0.5em} \text{$\forall$ $j \in \crb{1,...,N_k+1}$, $\forall$ $k \in \crb{1,...,K}$},
\end{equation}
where $c_i$ is the $i$-th SVIC for all $i = 1,...,n_c$ with $n_c$ being the total number of SVICs and $c_{i,\max/\min}$ is the corresponding upper and lower limits of the SVIC, respectively. Once the relative difference in Eq.~\eqref{eq:svic_rel_diff} is obtained at every collocation point plus the final non-collocated point, the structure detection algorithm for estimating the activation/deactivation points is executed as follows:

\begin{shadedframe}
\vspace{-10pt}
\begin{center}
 \shadowbox{\bf Detection Algorithm for State-Inequality Path Constrained Optimal Control Problems} 
\end{center}
\begin{enumerate}[{\bf Step 1:}]
\item Set detection tolerance for which renders a constraint active: $\epsilon$.
\item Check relative difference to identify activation/deactivation points: $\forall$ $ j \neq \{1\}$, $\forall$ $ k \neq \{1\}$
	\begin{enumerate}[{\bf (a):}]
	\item If $\delta c_i (\tau_j^{(k)}) \leq \epsilon$, and 
		\begin{enumerate}[{\bf (i):}]
		\item $\delta c_i (\tau_{j-1}^{(k)}) > \epsilon$, and $\delta c_i (\tau_{j+1}^{(k)}) \leq \epsilon$, then $\tau_j^{(k)} = \tau_s$ is deemed an \textit{activation} point. \vspace{0.5em}
		\item $\delta c_i (\tau_{j-1}^{(k)}) \leq \epsilon$, and $\delta c_i (\tau_{j+1}^{(k)}) > \epsilon$, then $\tau_j^{(k)} = \tau_s$ is deemed a  \textit{deactivation} point.				
		\end{enumerate}	 	
	\end{enumerate}
\item Compute bounds on estimated activation/deactivation points: $\left( \tau_s \right)^{\pm}$ 
\end{enumerate}
\end{shadedframe} 
The default detection tolerance that is used in the detection algorithm is chosen such that it follows a process analogous to selecting the NLP solver tolerance. Specifically, the NLP tolerance is typically chosen to be on the order of the square root of machine tolerance ($\sim \mathcal{O}(10^{-8})$), and so the default detection tolerance is chosen similarly to be on the order of the square root of the NLP solver tolerance ($\sim \mathcal{O}(10^{-4})$) given the solution on the initial mesh is typically governed by the accuracy at which the NLP solver can achieve.
Lastly, it is noted that the final step in the detection algorithm is performed because the detected A/D times are included as additional decision variables in the resulting nonlinear programming problem (NLP) and must be bounded to prevent overlapping or collapsing of domains within the partitioned problem. The bounds are obtained by
\begin{equation}\label{eq:ADbounds}
  \left.
    \begin{array}{lcl}
      \tau_s^-  &=& \tau_j^{(k)} + \nu \left( \tau_{j-1}^{(k)} - \tau_j^{(k)} \right), \\ \\
      \tau_s^+ &=&  \tau_j^{(k)} + \nu \left( \tau_{j+1}^{(k)} - \tau_j^{(k)} \right),
    \end{array} \right\} \left. \begin{array}{ll}  \\ s \in \{1,\ldots,S\}, \\ \,\end{array} \right.
\end{equation}
where $S = D - 1$ is the number of detected A/D times and $\nu > 0$ is a user defined parameter used to control the size of the allowable search space for the NLP solver to determine the optimal A/D times. The final step is to use the affine transformation given by Eq.~\eqref{eq:affine_map} to map the estimated A/D times back onto the original time domain $t \in [t_0,t_f]$. 
\subsubsection{Structure Decomposition}\label{sec:structure_decomp}
Assuming the structure detection method of Sect.~\ref{section:structure_detect} has identified activation and deactivation (A/D) times in the SVICs, the next step is to decompose the initial mesh into the multiple-domain structure. Specifically, the detected A/D times directly correlate to the number of domain interface variables, $t_s^{\crb{d}}$, $d \in \crb{1,...,D-1}$, to be solved for on subsequent mesh iterations, where the initial guess for each variable is the estimated
A/D times $\tau_s$, $s \in \crb{1,...,S}$ obtained from the method of Sect.~\ref{section:structure_detect}. The domain interface variables are included as additional decision variables in the NLP decision vector. Recall, the domain interface variables define the start of new domains, $\C{P}^{\crb{d}} = [t_s^{\crb{d}},t_s^{\crb{d-1}}]$, $d \in \crb{1,...,D-1}$, and are used to partition the time horizon $t \in [t_0,t_f]$ of the original optimal control problem into $D$ domains as described in Sect.~\ref{section:MDLGR}. \par
The next step is to enforce an additional constraint on the domain interface variables by including lower and upper bounds to prevent collapsing or overlapping of subdomains. The lower and upper bounds on each interval are determined by taking the bounds given by Eq.~\eqref{eq:ADbounds} and transforming them to the time interval $t \in [t_0,t_f]$ using the affine transformation of Eq.~\eqref{eq:affine_map}. Thus, the bounds for each domain interface variable $[\tau_s^-,\tau_s^+]$, $s \in \crb{1,...,S}$ are transformed to $[t_l^{\crb{d-1}},t_u^{\crb{d}}]$, $d \in \crb{1,...,D-1}$. \par
\begin{figure}
\begin{center}
	\includegraphics[scale=1]{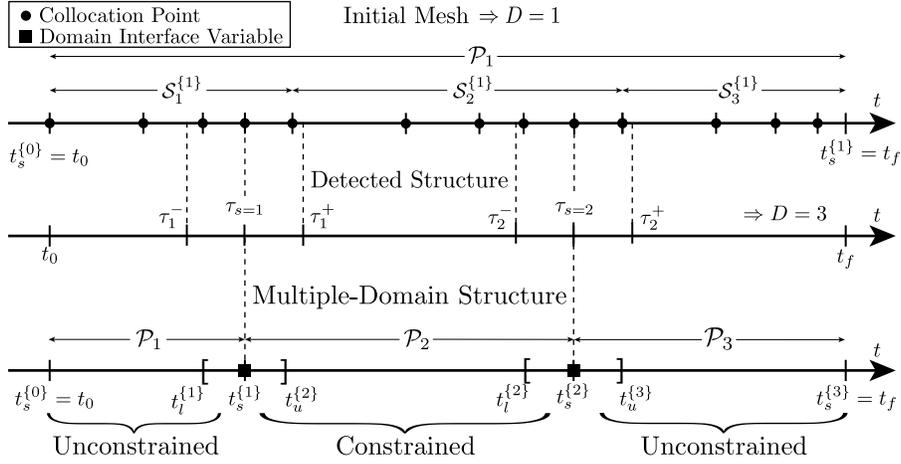}
	\caption{Schematic of process for decomposing the original optimal control problem (with $K=3$ mesh intervals) into $D$ domains where the $D-1$ domain interface variables are included as optimization variables to determine the optimal activation and deactivation in each SVIC.\label{fig:structureDecomp}}
\end{center}
\end{figure}

The overall approach for partitioning the time horizon of the original optimal control problem into multiple subdomains is shown in Fig.~\ref{fig:structureDecomp}, which provides an example of the initial domain consisting of $K=3$ mesh intervals. Although the domain interface variables are used to classify each subdomain, their precise locations are not fixed. The process of including the switch times as additional decision variables in the resulting NLP decision vector allows for their precise locations to be optimized. Recall, the method developed in this research assumes the SVICs are active over a non-zero interval of time, meaning no touch points exist. Thus, the SVIC will be active over a defined interval of time which will produce \textit{at least} one subdomain where the constraint(s) is active (as depicted in Fig.~\ref{fig:structureDecomp}). It is noted here that if the maximum relative error computed on the initial mesh ($D = 1$) does not meet the specified mesh tolerance, a mesh refinement method is used to obtained a suggested refined mesh. The provided refined mesh is then used when defining the mesh intervals on the partitioned problem.

\subsection{Domain Constraints and Refinement}\label{section:constraints_and_refinement}
Now that the time domain of the previous mesh has been partitioned into subdomains by the structure detection and decomposition method of Sect.~\ref{section:structure_detect_decomp}, the next step is to enforce additional necessary conditions in the subdomains where SVICs have been identified to be active. Once the additional necessary conditions described in Sect.~\ref{section:additional_nc} are properly enforced, the resulting nonlinear programming problem (NLP) is solved. Lastly, the maximum relative error and the maximum constraint violation is analyzed to determine if the domain refinement methods of Sect. \ref{section:domain_refinement_reg} or Sect.~\ref{section:domain_refinement_const} must be employed, respectively.
\subsubsection{Additional Necessary Conditions}\label{section:additional_nc}
As described in Sect.~\ref{section:svic}, the method developed in this work employs the method of indirect adjoining for handling the presence of state-variable inequality path constraints (SVICs). The method of indirect adjoining requires that the $q$-th time derivative of the SVIC be adjoined to the augmented Hamiltonian of the system (as shown in Eq.~\eqref{eq:adjoined_Hamiltonian}). In a direct method for discretizing an optimal control problem, the method of indirect adjoining is equivalent to enforcing the $q$-th time derivative of an active SVIC as an equality constraint throughout the constrained domain
\begin{equation}\label{eq:mixed_equality_constraint}
\dfrac{\text{d}^{(q)}c}{\text{d}t^{(q)}} = c^{(q)}\left( \m{y}\left( t\hspace{0.035cm} \right),\m{u} \left(t\hspace{0.035cm}\right), t\hspace{0.035cm} \right) = 0, \hspace{0.5em} t \in [t_a,t_d]
\end{equation}
where $t_a$ and $t_d$ are the activation and deactivation times of the SVIC, respectively. Additionally, given that the $q$-th time derivative of the SVIC is being used to derive the optimal control, the path entering the constrained arc must satisfy the necessary tangency constraints (given by Eq.~\eqref{eq:tangencyconstraints}) \cite{Bryson1963}. In summary, the additional necessary conditions that are added to the NLP presented in Section~\ref{section:NLP} is (i) the mixed state-control equality constraint on each collocation point within a constrained domain and (ii) the necessary tangency constraints on the first collocation point of the first mesh interval within a constrained domain. \par
\begin{figure}[htp]
\begin{center}
	\includegraphics[scale=0.7]{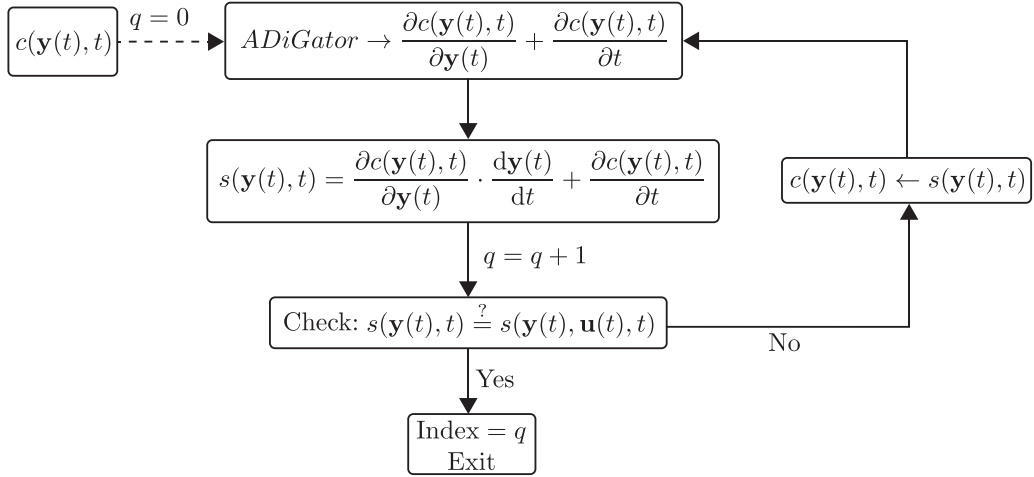}
	\caption{Schematic of algorithm for obtaining higher-order time derivatives of any active state-variable inequality path constraint.\label{fig:IndexID}}
\end{center}
\end{figure}

 
In this work it is assumed that the functions defining the SVICs are continuous and differentiable, and although the method developed in this paper is not restricted on how the derivatives of the SVICs are obtained, the process in which the time derivatives are computed in this work is described as follows. To avoid having to derive higher-order time derivatives analytically, algorithmic differentiation techniques and the properties of the chain rule are exploited to automate the process. Figure~\ref{fig:IndexID} provides a schematic of how the higher-order time derivatives of the active SVICs are obtained. The algorithmic differentiation software $ADiGator$ \cite{weinstein2017algorithm} is used to obtain the partial derivatives of  the active SVIC with respect to both the state and time. To complete the full time derivative, a stand-alone function $s(\m{y}(t),t)$ is written by applying the chain rule and properly substituting the dynamic constraints $\m{f}(\m{y}(t),\m{u}(t),t)$ for $\text{d}\m{y}(t)/\text{d}t$. The algorithm then checks if the stand-alone function is now a function of any control variables. If the stand-alone function is a function of any of the control variables then the algorithm terminates. Otherwise, the stand-alone function re-enters the algorithm loop and is treated as if it were the original SVIC (given it is not yet a function of any of the control variables). Upon completion, the derivative count $q$ defines the order of the SVIC and is used when defining the additional necessary conditions given by Eqs.~\eqref{eq:tangencyconstraints} and \eqref{eq:mixed_equality_constraint}. Note, the process shown in Fig.~\ref{fig:IndexID} only needs to be executed once and can be done before solving the problem.

\subsubsection{Regular Domain Refinement}\label{section:domain_refinement_reg}
Suppose now that the solution $\m{Y}^{(M)}$ obtained after performing structure detection and enforcing the additional necessary conditions satisfies the condition $\max \{ \Delta c_i^{(M)} \} \leq \epsilon_{c}$, that is, the maximum constraint violations are within the user specified constraint satisfaction tolerance. Although the constraints have met the specified tolerance, certain domains may still require mesh refinement. That is, the maximum relative error $e^{(M)}_{\max}$ from the solution $\m{Y}^{(M)}$ exceeds the specified mesh error tolerance. To date, several mesh refinement methods have been developed \cite{Darby:2011a,Darby:2011b,Patterson:2015,Liu:2017,Liu2018}. For the problems solved in this work, the domains which did not meet the mesh tolerance were refined using the mesh refinement method developed in Ref. \cite{Liu2018}. The process for which Ref. \cite{Liu2018} computes the maximum relative error follows the same process as shown in Sect. \ref{section:solution_error}. The reader is referred to Ref. \cite{Liu2018} for a detailed explanation as to how the mesh refinement method operates.

\subsubsection{Constrained Domain Refinement}\label{section:domain_refinement_const}
Suppose $\m{Y}^{(M)}$ is the solution after solving the reformulated NLP obtained from performing the structure detection and decomposition method of Sect. \ref{section:structure_detect_decomp} and enforcing the additional necessary conditions as given in Sect. \ref{section:additional_nc}. Suppose now that the computed maximum constraint violation 
\begin{equation}
\max \{ \Delta c_i^{(M)} \} = \max \{ c_i(\m{Y}^{(M)}) - c_{i,\max} \} 
\end{equation}
on the current mesh $M$ is larger than the user specified constraint satisfaction tolerance $\epsilon_{c}$, that is, the condition $\max \{ \Delta c_i^{(M)} \} \geq \epsilon_{c}$ holds. Suppose further that the maximum relative error on the current mesh $e^{(M)}_{\max}$ does not satisfy the specified mesh tolerance $\epsilon_{mesh}$, that is, the condition $e^{(M)}_{\max} > \epsilon_{mesh}$ holds. It may be the case that the reformulated NLP is ill-conditioned and/or the grid on which the problem was previously solved is too coarse. Ill-conditioning can result from enforcing the additional necessary conditions at incorrect locations throughout the trajectory, which in this work would imply the detected structure could be inaccurate.  To account for this possibility, the next step taken is to re-perform structure detection using the solution $\m{Y}^{(M)}$. Additionally, given that the mesh error tolerance has not been met, a mesh refinement method (for example, Ref. \cite{Liu2018}) is used to refine the previous mesh. The problem is now reformulated using the new detection results while the mesh intervals within each domain is specified in accordance with the newly refined mesh provided by the chosen mesh refinement method.

\subsection{Procedure for State-Path Constrained Optimal Control Problems}\label{section:procedure}

An overview of the method developed in this work is provided below. The mesh iteration is denoted by $M$ and is incremented by one upon each loop of the method. The method terminates when either a maximum prescribed mesh count, $M_{\max}$, is reached or when both the mesh error tolerance, $\epsilon_{mesh}$, is satisfied on every mesh interval and the maximum constraint violation for each active SVIC, $\max\{\Delta c_{i}^{(M)}\}$ $i = 1,...,n_c$, is less than the user specified constraint satisfaction tolerance ($\epsilon_{c}$).
\begin{shadedframe}
\vspace{-10pt}
\begin{center}
 \shadowbox{\bf Method for Solving State-Path Constrained Optimal Control Problems} 
\end{center}
\begin{enumerate}[{\bf Step 1:}]
\setlength\itemsep{1em}
\item Set $M = 0$ and specify initial mesh. All mesh intervals form a single domain.
\item Solve the NLP resulting from Radau collocation presented in Section \ref{section:NLP} on mesh $M = 0$. \\
\hdashrule[0.5ex]{15cm}{0.5pt}{1pt}
\vspace{-1em}
\item Employ the structure detection and decomposition method of Section \ref{section:structure_detect_decomp}. \vspace{0.5em}
	\begin{enumerate}[{\bf (a):}]
	\setlength\itemsep{1em}
		\item Determine activation/deactivation times in each constraint using the method of Section \ref{section:structure_detect}. 
		\item Partition the time horizon into domains using the method of Section \ref{sec:structure_decomp}.
		\item Classify the domains as either constrained or unconstrained.
		\item Enforce additional necessary conditions as described in  Section \ref{section:additional_nc}.
	\end{enumerate}	
\item Increment $M \rightarrow M+1$ and solve the reformulated NLP.
\item Compute maximum constraint violations, $\max\{\Delta c_{i}^{(M)}\}$, and maximum relative error, $e_{\max}^{(M)}$, on current mesh.
\item If ($\max\{\Delta c_{i}^{(M)}\} \leq \epsilon_{c}$ and $e_{\max}^{(M)} \leq \epsilon_{mesh}$) or ($M > M_{\max}$), then quit. Otherwise, \vspace{0.5em}
\begin{enumerate}[{\bf (a):}]
\setlength\itemsep{1em}
\item If $\max\{\Delta c_{i}^{(M)}\} \geq \epsilon_{c}$,
\begin{enumerate}[{\bf (i):}]
\item Perform domain refinement according to Section \ref{section:domain_refinement_const}
\item Return to \textbf{Step 3}.
\end{enumerate}
\item If $\max\{\Delta c_{i}^{(M)}\} \leq \epsilon_{c}$, and $e_{\max}^{(M)} > \epsilon_{mesh}$,
\begin{enumerate}[{\bf (i):}]
\item Perform domain refinement according to Section \ref{section:domain_refinement_reg}
\item Return to \textbf{Step 4}.
\end{enumerate}
\end{enumerate}
 
\end{enumerate}
\end{shadedframe} 

%%%%%%%%%%%%%%%%%%%%%%%%%%%%%%% SECTION 6 %%%%%%%%%%%%%%%%%%%%%%%%%%%%%%%%%%%%%%%%%%%%%%%%%%%%%%%
\section{Examples}\label{section:examples}
In this section, the aforementioned SPOC method described in Sect. \ref{section:method} is demonstrated on two examples. The first example consists of a second-order state-variable inequality path constraint with a known analytical solution which provides a baseline to analyze the accuracy of the proposed method. The second example involves a problem with a first-order state-variable inequality path constraint where an analytical solution is not available and demonstrates the methods ability to improve constraint satisfaction when analyzing a nonlinear dynamical system. In both examples, the performance of the SPOC method developed in this paper will be compared against the mesh refinement method developed in Ref. \cite{Liu2018}. Given the mesh refinement method developed in Ref. \cite{Liu2018} is implemented within the SPOC method, comparisons were chosen to be made against this specific mesh refinement method. Additionally, while comparisons could be made against other mesh refinement methods (such as the methods of Ref. \cite{Darby:2011a} and Ref. \cite{Patterson:2015}), it has been found that the method of Ref. \cite{Liu2018} typically outperforms these previously developed mesh refinement methods. Lastly, it is noted that the method of Ref. \cite{Liu2018} will be referred to from this point forth as the $hp$-LGR method.\par
All results to follow will include three solutions. The first being the solution obtained in the first two steps of the procedure given in Sect. \ref{section:procedure}, referred to from this point forth as the initial mesh. The second solution shown is obtained using the $hp$-LGR method and the third solution is obtained using the SPOC method. The design parameters used within the method were chosen for both examples to be: $\epsilon = 1\times10^{-4}$ and $\nu = 5$. The specified mesh tolerance for both the $hp$-LGR method and the SPOC method was chosen to be $\epsilon_{mesh} = 1 \times 10^{-6}$ for the first example, while for the second example, the mesh tolerance for both methods was set to $\epsilon_{mesh} = 1 \times 10^{-8}$. All results are obtained, other than those obtained using the SPOC method, using the \textsf{MATLAB}$^\textsf{\textregistered}$ optimal control software $\mathbb{GPOPS-II}$ as described in Ref.~\cite{Patterson2014}. The resulting nonlinear programming problem (NLP) was solved using the NLP solver IPOPT~\cite{Biegler2008} in full-Newton mode with an NLP tolerance and NLP constraint violation tolerance both set to $\epsilon_{NLP} = 1 \times 10^{-8}$. The user specified constraint satisfaction tolerance for both examples was set equal to the NLP tolerance, $\epsilon_c = \epsilon_{NLP}$. The first and second derivatives supplied to IPOPT were obtained using the built-in sparse finite central difference method in $\mathbb{GPOPS-II}$ which uses the method of Ref. \cite{patterson2012exploiting}. In both examples, the initial mesh consisted of ten mesh intervals with four collocation points within each mesh interval. The initial guess for the initial mesh consisted of a straight line connecting boundary conditions for states with boundary conditions specified at both the beginning and end of the original time domain, while a constant value was supplied for states with a single boundary condition specified at either end of the original time domain. Finally, all numerical computations were performed on a 2.9 GHz Intel Core i9 MacBook Pro running macOS Monterey Version 12.5 with 32 GB 2400 MHz DDR4 of RAM, using \textsf{MATLAB}$^\textsf{\textregistered}$ version R2021a (build 9.10.0.1669831).  
\subsection{Example 1: Bryson-Denham Problem}
\begin{figure}[htp]
\begin{center}
	\includegraphics[scale=0.4]{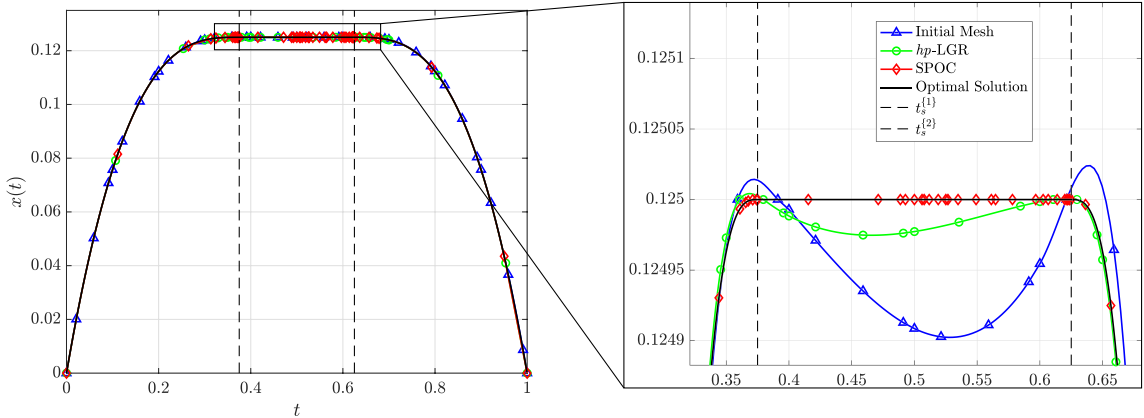}
	\caption{Position, $x(t)$, comparison between the solution obtained on the initial mesh, with the $hp$-LGR method of Ref. \cite{Liu2018}, and the SPOC method for the Bryson-Denham problem.\label{fig:bd_position}}
\end{center}
\end{figure}

Consider the following optimal control problem obtained from Ref. \cite{Bryson1963}, often referred to as the Bryson-Denham problem, which includes a second-order state-variable inequality path constraint. Minimize the objective functional
\begin{equation}
\C{J} = \dfrac{1}{2} \int_{0}^{1} u^2(t) \text{ d}t,
\end{equation}
subject to the dynamic constraints
\begin{equation}
\begin{array}{lcl}
\dot{x}(t) &=& v(t),\\
\dot{v}(t) &=& u(t),
\end{array}
\end{equation}
the boundary conditions
\begin{equation}
\begin{array}{lclcl}
x(0) &=& x(1) &=& 0, \\
v(0) &=& -v(1) &=& 1,
\end{array}
\end{equation}
and the state-variable inequality path constraint (SVIC)
\begin{equation}
x(t) \leq L,
\end{equation}
where $L$ is the upper limit on the SVIC. Though the problem formulation may appear rudimentary, it is challenging to obtain an accurate numerical approximation of the solution. For completeness, the analytical solution for the control, position, and objective value to the above optimal control problem for $0 \leq L \leq 1/6$ is given by
\begin{align}
u^*(t) &= 
\begin{cases} 
-\dfrac{2}{3L} \left(1 - \dfrac{t}{3L} \right), & 0 \leq t \leq 3L \vspace{0.3em}\\[2pt]
0, & 3L \leq t \leq 1 - 3L \vspace{0.3em}\\[2pt]
-\dfrac{2}{3L} \left(1 - \dfrac{1-t}{3L} \right), & 1-3L \leq t \leq 1 \\
\end{cases} \\[5pt] 
x^*(t) &= 
\begin{cases}
L \left[ 1 - \left(1 - \dfrac{t}{3L}\right)^3 \right], & 0 \leq t \leq 3L \vspace{0.3em}\\[2pt]
L, & 3L \leq t \leq 1 - 3L \vspace{0.3em}\\[2pt]
L \left[ 1 - \left(1 - \dfrac{1-t}{3L}\right)^3 \right], & 1-3L \leq t \leq 1\\
\end{cases} \\[5pt]
\C{J}^* &= \dfrac{4}{9L}.
\end{align} 
For the results to follow, the upper limit on the SVIC was chosen to be $L = 1/8 \Rightarrow \C{J}^* = 32/9$. Figure~\ref{fig:bd_position} shows the computed position, $x(t)$, for three solutions: the first obtained on the initial mesh, the second computed using the $hp$-LGR method, and the third solution computed using the SPOC method. Upon zooming in on the portion of the trajectory where the SVIC is known to be active, it is observed that both solutions obtained on the initial mesh and using the $hp$-LGR method do not remain on the constrained arc. Conversely, the position computed using the SPOC method is able to accurately approximate the active constraint arc. 
\begin{figure}[htp]
\begin{center}
	\includegraphics[scale=0.4]{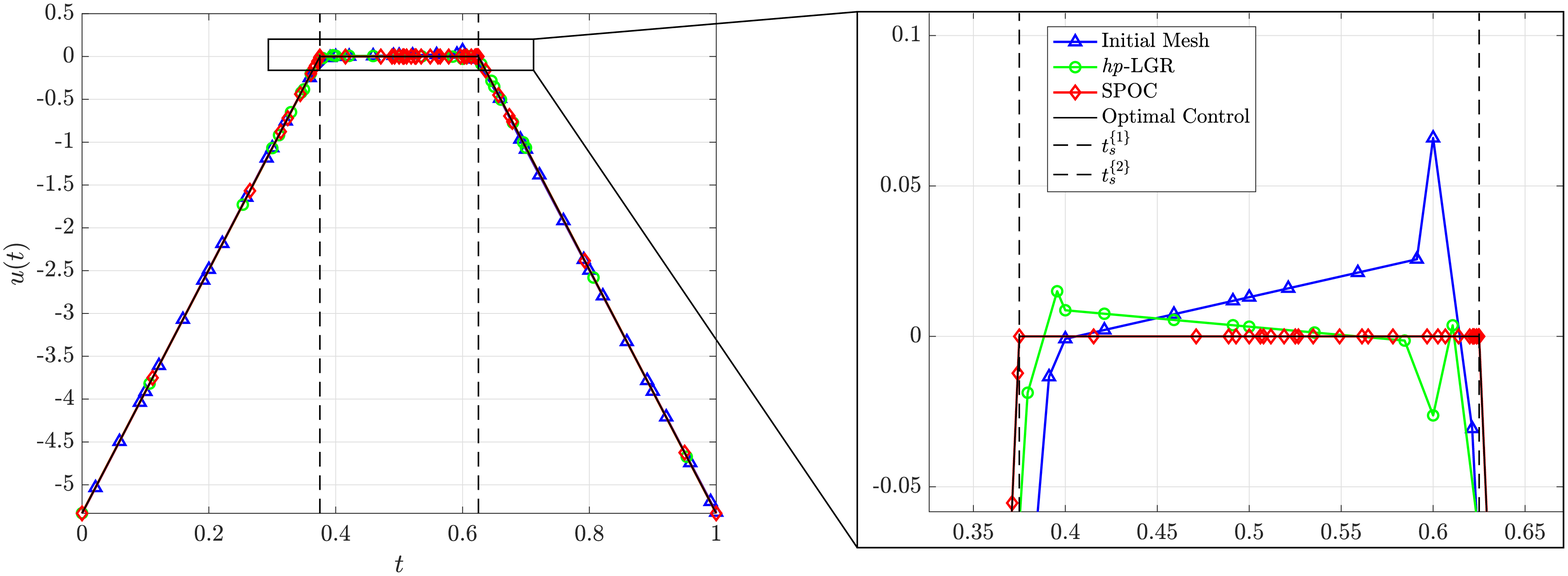}
	\caption{Control, $u(t)$, comparison between the solution obtained on the initial mesh, with the $hp$-LGR method of Ref. \cite{Liu2018}, and the SPOC method for the Bryson-Denham problem.\label{fig:bd_control}}
\end{center}
\end{figure}
A similar comparison is given in Fig.~\ref{fig:bd_control} for the time-history of the computed optimal control, $u(t)$. Specifically, upon zooming in on the region where the constraint is known to be active, Fig.~\ref{fig:bd_control} shows that the optimal control computed on the initial mesh and with the $hp$-LGR method exhibits "chattering" behavior near the activation and deactivation of the SVIC. Conversely, the optimal control computed by the SPOC method removes any chattering behavior and closely approximates the known optimal control along the active constraint arc.\par

Using a detection tolerance of $\epsilon =1 \times 10^{-4}$, and a value of $\nu = 5$, Table~\ref{table:estimatedtimes} provides the estimated activation and deactivation (A/D) times ($\hat{t}_s$) using the detection algorithm presented in Sect. \ref{section:structure_detect}. Also provided in Table~\ref{table:estimatedtimes} are the lower and upper bounds ($\hat{t}_s^-$, $\hat{t}_s^+$) on the detected times, respectively, the optimized A/D times ($t_s^*$) returned by the NLP solver, and the known analytic switch times. By comparing the estimated and optimized switch times with the analytical switch times, it is observed that the NLP solver is capable of selecting switch times which are in close agreement to the known analytical A/D times. This result shows that, given the provided search space for the A/D times contains the optimal switch time, the NLP solver is capable of optimizing the location of the A/D times such that it is in close agreement with the analytical solution. Note, the dashed lines shown in both Figs. \ref{fig:bd_position} and \ref{fig:bd_control} represent the activation and deactivation times returned by the NLP. 
\begin{table}[htbp]
\centering
\caption{Detection performance: Estimated A/D times, $\hat{t}_s$, lower and upper bounds $\hat{t}_s^-,\hat{t}_s^+$, optimized A/D times $t_s^*$, and the analytical A/D times}
\begin{tabular}{ l c c c c c} 
 \hline \hline
 \rule{0pt}{3ex}Switch time & $\hat{t}_s$ & $\hat{t}_s^{-}$ & $\hat{t}_s^{+}$ & $t_s^*$ & Analytic\\[3pt] \hline 
 \rule{0pt}{3ex}Activation & 0.35905 & 0.16996 & 0.51949 & 0.37498 & 0.3750\\[3pt]
 \rule{0pt}{3ex}Deactivation & 0.65905 & 0.46996 &  0.81949 &  0.62501 & 0.6250\\[3pt] 
 \hline \hline
\end{tabular}
\label{table:estimatedtimes}
\end{table}
Table~\ref{table:relerror} provides the CPU times and computed relative errors with respect to the known analytical solution in the objective value, $\delta \C{J}$, the activation time, $\delta  t_s^{[1]}$, and the deactivation time, $\delta  t_s^{[2]}$, for each solution shown in Figs.~\ref{fig:bd_position} and \ref{fig:bd_control}. 
\begin{table}[htp]
\centering
\caption{Comparison of the CPU time and the relative error of the computed cost, $\delta \C{J}$, the activation time, $\delta t_s^{\{1\}}$, and the deactivation time, $\delta t_s^{\{2\}}$, with respect to the analytical solution for the initial mesh, $hp$-LGR method, the UTM method \cite{Mall2020}, and the SPOC method}
\begin{tabular}{ l l c c c c c} 
 \hline \hline
 \rule{0pt}{3ex} & \multicolumn{1}{c}{$\delta \C{J}$} & $\delta t_s^{\{1\}}$ $[s]$ & $\delta t_s^{\{2\}}$ $[s]$ & CPU Time $[s]$\\[3pt] \hline  
 \rule{0pt}{3ex}Initial Mesh & 4.19$\times 10^{-5}$ & 3.63$\times 10^{-2}$ & 1.14$\times 10^{-2}$ & 0.32\\[3pt] 
 \rule{0pt}{3ex}$hp$-LGR & 9.72$\times 10^{-6}$ & 2.88$\times 10^{-2}$ &  1.79$\times 10^{-2}$ & 0.40 \\[3pt] 
 \rule{0pt}{3ex}UTM \cite{Mall2020} & 3.78$\times 10^{-4}$ & -- &  -- & -- \\[3pt] 
 \rule{0pt}{3ex}SPOC & 3.02$\times 10^{-13}$ & 2.03$\times 10^{-5}$ &  8.85$\times 10^{-6}$ & 2.36\\[3pt] 
 \hline \hline
\end{tabular}
\label{table:relerror}
\end{table}
Additionally, Table~\ref{table:relerror} includes a comparison against the Uniform Trigonometrization Method (UTM) (developed in Ref. \cite{Mall2020}) which follows the procedures of a penalty method within an indirect formalism. It is found that the SPOC method computes an objective that is eight orders of magnitude more accurate than the objective obtained using the $hp$-LGR method, and up to nine orders of magnitude of improvement with respect to the objective given in Ref. \cite{Mall2020}. Additionally, an improvement of up to four orders of magnitude in the computed A/D times is obtained when using the SPOC method.

\subsection{Example 2: Constrained Entry Vehicle Crossrange Maximization}
Consider the following constrained reusable launch vehicle reentry (RLVE) optimal control problem, obtained from Ref.~\cite{Betts2010}, of maximizing the achieved crossrange during atmospheric entry while being subject to a stagnation point heating rate constraint. Minimize the cost functional
\begin{equation}
\C{J} = -\phi(t_f)
\end{equation}
subject to the dynamic constraints
\begin{subequations}\label{eq:EOM}
\begin{align}
\dot{r} &= v \sin \gamma, \label{eq:rdot}\\[5pt]
\dot{\theta} &= \dfrac{v \cos \gamma \sin \psi}{r \cos \phi}, \label{eq:thetadot}\\[5pt]
\dot{\phi} &= \dfrac{v \cos \gamma \cos \psi}{r}, \label{eq:phidot}\\[5pt]
\dot{v} &= -D -g\sin \gamma, \label{eq:vdot}\\[5pt]
\dot{\gamma} &= \dfrac{L \cos \sigma}{v} + \cos \gamma \left( \dfrac{v}{r} - \dfrac{g}{v}  \right), \label{eq:gammadot}\\[5pt]
\dot{\psi} &=  \dfrac{L \sin \sigma}{v \cos \gamma} + \dfrac{v}{r} \cos \gamma \sin \psi \tan \phi,\label{eq:psidot}
\end{align}
\end{subequations}
the stagnation point heating rate constraint \cite{Detra1957}
\begin{equation}\label{eq:heatrate_eq1}
\dot{Q} = \hat{\dot{Q}} \left(\rho / \rho_0 \right)^{1/2} \left(v/v_c \right)^{3.15} \leq \dot{Q}_{\max},
\end{equation}
and the boundary conditions
\begin{equation}
\begin{array}{lclclcl}
h(t_0) &=& 79.248 \text{ km} &,& h(t_f) &=& 24.384 \text{ km}, \\[5pt]
\theta(t_0) &=& 0 \text{ deg} &,& \theta(t_f) &=& \text{Free}, \\[5pt]
\phi(t_0) &=& 0 \text{ deg} &,& \phi(t_f) &=& \text{Free}, \\[5pt]
v(t_0) &=& 7802.88 \text{ m/s} &,& v(t_f) &=& 762.0 \text{ m/s}, \\[5pt]
\gamma(t_0) &=& -1 \text{ deg} &,& \gamma(t_f) &=& -5 \text{ deg}, \\[5pt]
\psi(t_0) &=& 90 \text{ deg} &,& \psi(t_f) &=& -90 \text{ deg},
\end{array}
\end{equation}
where the terminal time $t_f$ is free, and the two control inputs are the vehicles angle of attack command, $\alpha(t)$ and the bank angle command, $\sigma(t)$. The lift and drag specific forces are modeled, respectively, as
\begin{equation}
  \begin{array}{lcl}
    D &=& q S C_D / m, \\
    L &=& q S C_L / m, 
  \end{array}
\end{equation}
with $q = \rho v^2 / 2$ being the dynamic pressure
\begin{equation}\label{eq:density}
  \rho = \rho_0 \exp \left(-\beta h \right)
\end{equation}
is the atmospheric density, and $h = r - R_e$ is the altitude of the vehicle. The lift and drag aerodynamic coefficients are modeled, respectively, as \cite{Betts2010}
\begin{equation}
\begin{array}{cll}
C_L &=& C_{L0} + C_{L1}\alpha, \\[5pt]
C_D &=& C_{D0} + C_{D1}\alpha + C_{D2}\alpha^2,
\end{array}
\end{equation}
The aerodynamic coefficients and physical constants used in this model are provided in Table~\ref{table:parameters}, noting that Ref.~\cite{Betts2010} uses English units whereas SI units were used to obtain the results presented in this work. Given that $\dot{Q}$ defined in Eq.~\eqref{eq:heatrate_eq1} is always positive during atmospheric flight, it can be transformed monotonically via the natural logarithm as
\begin{equation}\label{eq:logQdot}
\ell_{\dot{Q}} \equiv \log \dot{Q} =  \log  \hat{\dot{Q}}  - \beta h/2 + 3.15 \log(v/v_c),
\end{equation}
where the expression for the atmospheric density given in Eq.~\eqref{eq:density} was used to obtain the result shown above in Eq.~\eqref{eq:logQdot}.
Then the stagnation point heating rate constraint can be rewritten as
\begin{equation}\label{eq:heatrate_eq2}
\ell_{\dot{Q}} \leq \log \dot{Q}_{\max}.
\end{equation}
It is noted in Ref.~\cite{darby2011minimum} that Eq.~\eqref{eq:heatrate_eq2} behaves better computationally than the form given in Eq.~\eqref{eq:heatrate_eq1}.\par
\begin{figure}[htp]
\begin{center}
	\includegraphics[scale=0.35]{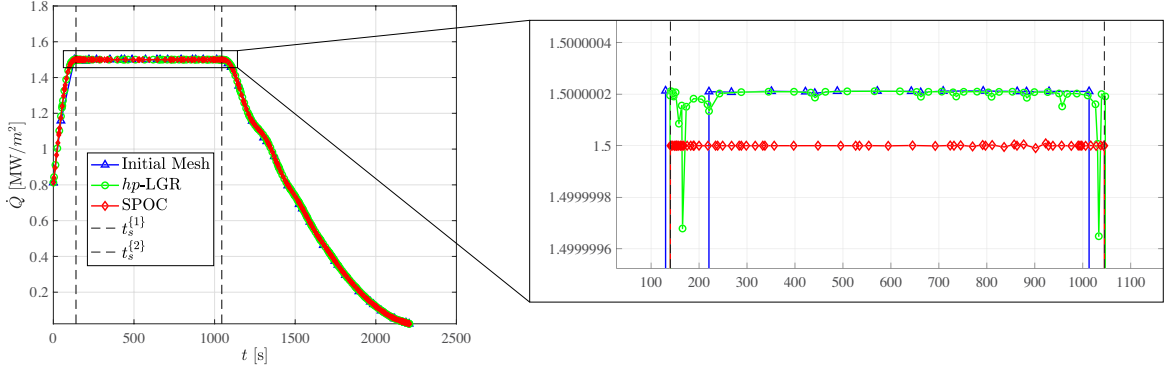}
	\caption{Heating rate profile, $\dot{Q}(t)$, comparison between the solution obtained on the initial mesh, with the $hp$-LGR method of Ref. \cite{Liu2018}, and the SPOC method for the constrained RLVE problem.\label{fig:qdot}}
\end{center}
\end{figure}

\begin{table}[h]
\centering
\caption{Physical and aerodynamic constants}
\begin{tabular}{ c c l } 
 \hline \hline 
 \rule{0pt}{2.5ex} Quantity & Value & Unit \\ \hline  
 \rule{0pt}{2.5ex} $R_e$ & 6371.20392 & km \\[3pt] 
 \rule{0pt}{2ex}  $1/\beta$ & 7254.24 & m  \\[3pt]
 \rule{0pt}{2ex}  $\rho_0$ & 1.22557083 & kg/m$^3$  \\[3pt]
 \rule{0pt}{2ex}  $\mu$ & 3.98603195 & m$^3$/s$^2$  \\[3pt]
 \rule{0pt}{2ex}  $m$ & 92079.2526 & kg  \\[3pt]
 \rule{0pt}{2ex}  $S$ & 249.909178 & m$^2$  \\
 \rule{0pt}{2ex}  $\hat{\dot{Q}}$ & 199.87$\times$10$^6$ & W/m$^2$  \\[3pt]
 \rule{0pt}{2ex}  $C_{L0}$ & $-0.2070$ & --  \\[3pt]
 \rule{0pt}{2ex}  $C_{L1}$ & 1.6756 & 1/rad  \\[3pt]
 \rule{0pt}{2ex}  $C_{D0}$ & $0.0785$ & --  \\[3pt]
 \rule{0pt}{2ex}  $C_{D1}$ & $-0.3529$ & 1/rad \\[3pt]
 \rule{0pt}{2ex}  $C_{D2}$ & $2.0400$ & 1/rad$^2$ \\[3pt]
  \hline \hline
\end{tabular}
\label{table:parameters}
\end{table}
For the results presented in this work, the maximum allowable stagnation point heating rate constraint was chosen to be $\dot{Q}_{\max} = 1.5 $ MW/m$^2$. Figure~\ref{fig:qdot} provides a comparison of the computed stagnation point heating rate profiles obtained (i) on the initial mesh (ii) using the $hp$-LGR method and (iii) using the SPOC method. Zooming in on the portion where the heating rate constraint looks to be active, a couple observations can be made. The first observation is that both heating rate profiles computed using the $hp$-LGR method and on the initial mesh retain a value that remains above the prescribed constraint limit. Furthermore, both heating rate profiles contain some level of chattering around the activation and deactivation points. Conversely, the stagnation point heating rate profile obtained using the SPOC method is able to remove any chattering behavior around the A/D times and retains a value that remains in close agreement with the specified upper limit of $\dot{Q}_{\max} = 1.5 $ MW/m$^2$. 

\begin{figure}
\centering
\vspace*{0.25cm}
\subfloat[\label{fig:mesh_error}]{\includegraphics[scale=0.4]{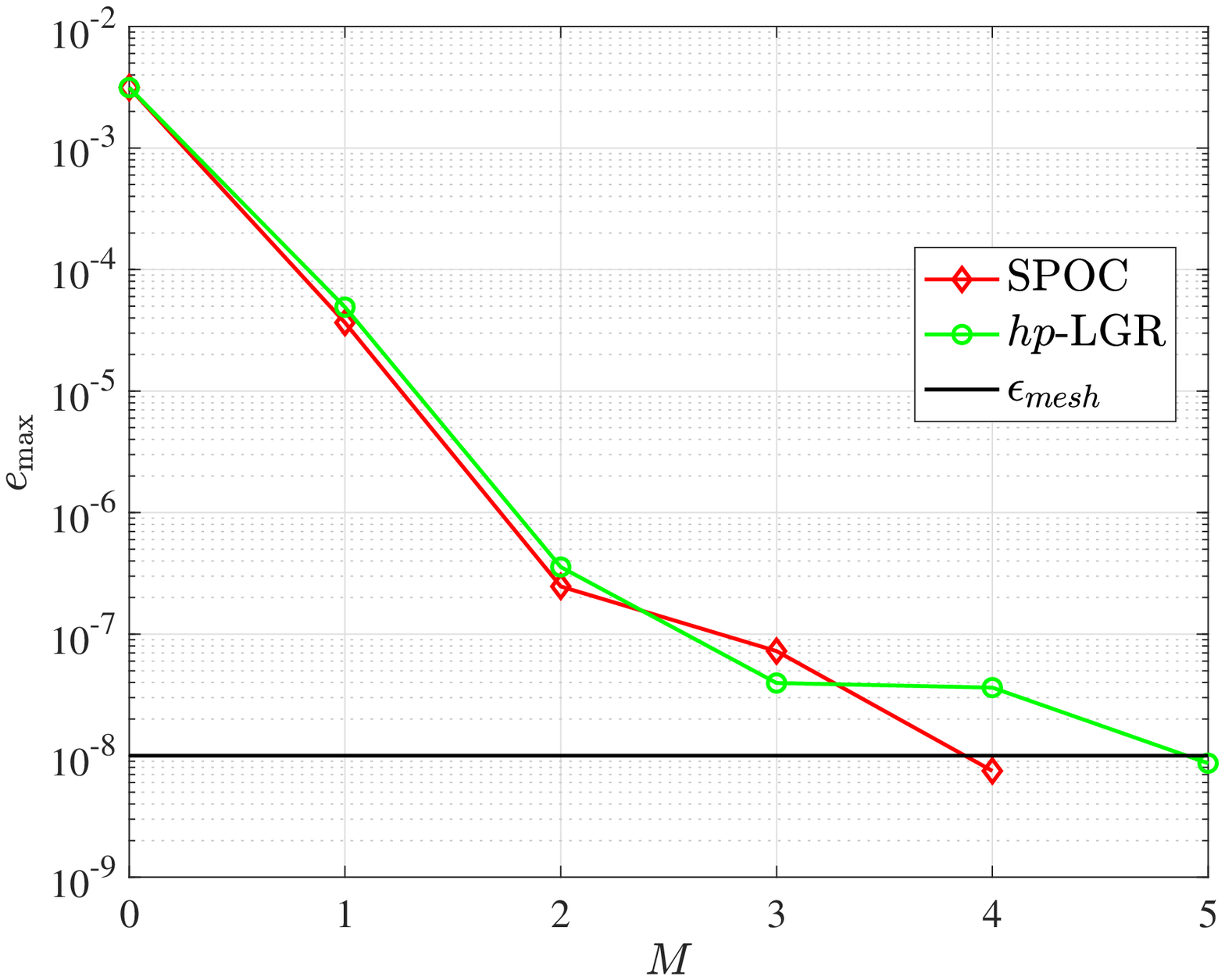}}~~\subfloat[\label{fig:constr_viol}]{\includegraphics[scale=0.4]{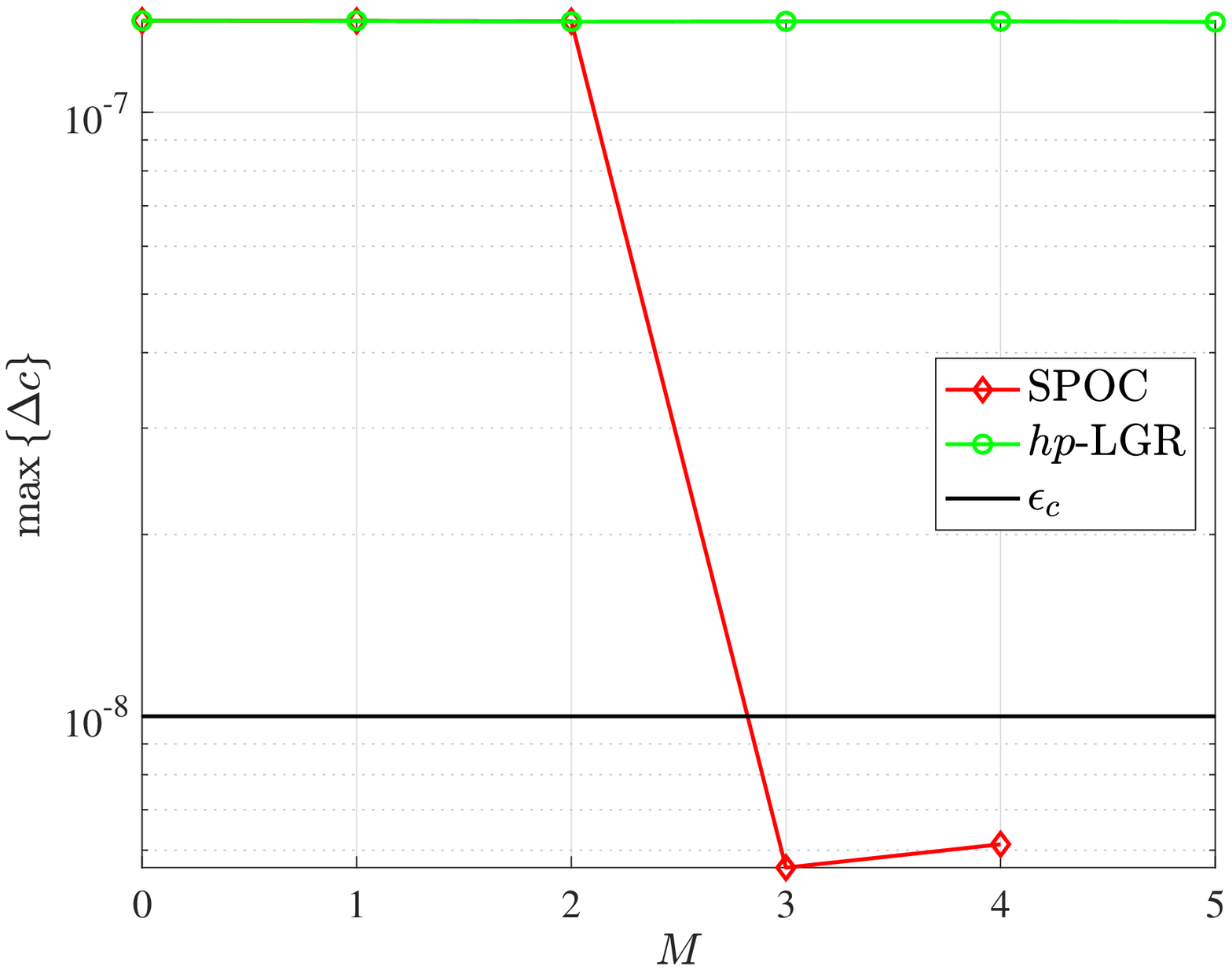}} \\ 
\caption{Maximum relative error, $e_{\max}$ and maximum constraint violations, $\max{\{ \Delta c\}}$ versus mesh iteration, $M$ for the $hp$-LGR and SPOC methods, respectively.} 
\label{fig:error_vs_mesh}
\end{figure}

Figure~\ref{fig:mesh_error} compares the computed maximum relative error (given by Eq.~\eqref{eq:max_mesh_err}) at each mesh iteration for both the $hp$-LGR method and the SPOC method. It is found that the SPOC method takes fewer mesh iterations to meet the specified error tolerance and on average obtains a maximum relative error that is smaller than the maximum relative error obtained with the $hp$-LGR method at each mesh iteration. Likewise, Fig.~\ref{fig:constr_viol} compares the computed maximum constraint violation during each mesh iteration for both the $hp$-LGR method and the SPOC method. It is observed that the SPOC method is capable of reducing the computed maximum constraint violation below the specified NLP tolerance while the maximum constraint violation computed using the the $hp$-LGR method remains approximately an order of magnitude above the specified NLP tolerance during each mesh iteration.

Table~\ref{table:performance_comp} provides a comparison of the computed objective values, maximum constraint violations, and CPU times for the solutions obtained (i) on the initial mesh (ii) using the $hp$-LGR method and (iii) using the SPOC method. First, it can be seen that the SPOC method computes an objective value that is marginally smaller than the objective obtained using the $hp$-LGR method. The decrease in cost is likely due to the fact that the solution obtained using the SPOC method satisfies the constraint to a tighter tolerance, and does not continuously violate the constraint throughout the constraint arc. As depicted in Fig.~\ref{fig:constr_viol}, Table~\ref{table:performance_comp} re-emphasizes that the SPOC method is capable of reducing the maximum constraint violation as compared with the initial mesh and the $hp$-LGR method. Although the SPOC method converged to a solution in one less mesh iteration than the $hp$-LGR method, Table~\ref{table:performance_comp} shows that the SPOC method took an order of magnitude longer to compute a solution. A large percentage of the CPU time was spent within the NLP solver, and may be due to the inclusion of additional necessary conditions and decision variables.

\begin{table}[htp]
\centering
\caption{Computed cost, maximum constraint violations, and CPU time for the initial mesh, $hp$-LGR method, and the SPOC method}
\begin{tabular}{ l c c c c c} 
 \hline \hline 
 \rule{0pt}{3ex} & \multicolumn{1}{c}{$\C{J} = \phi(t_f) [$deg$]$} & $\max \{ \dot{Q}(t) - \dot{Q}_{\max} \}$ $[$MW/m$^2]$ & CPU Time $[s]$\\[3pt] \hline   
 \rule{0pt}{3ex}Initial Mesh & 33.4325507 & 1.418$\times 10^{-7}$ &  1.34\\[3pt] 
 \rule{0pt}{3ex}$hp$-LGR & 33.4465173 & 1.411$\times 10^{-7}$ &   6.40 \\[3pt] 
 \rule{0pt}{3ex}SPOC & 33.4465161 & 6.137$\times 10^{-9}$ &  50.9\\[3pt] 
 \hline \hline
\end{tabular}
\label{table:performance_comp}
\end{table}

It is emphasized that the structure detection algorithm of Sect. \ref{section:structure_detect} is only performed when the maximum constraint violation is above the user specified constraint satisfaction tolerance $\epsilon_{c}$, and for this example, it was executed three times before the maximum constraint violation was reduced below $\epsilon_{c}$. A list of estimated and optimized A/D times obtained during the SPOC mesh iterations for which structure detection was employed is provided in Table~\ref{table:detection_performance}. Note, estimated A/D times serve as initial guesses for the optimized A/D times which are returned by the NLP solver after solving the partitioned problem. Observing the trend in the relative error, it can be seen that the estimated and optimized A/D times converge on one another as the mesh iteration increases. This is likely due to the fact that the solution for which detection is being performed on is (i) on a denser mesh and (ii) approaching the correct constraint structure. The ability to recursively perform structure detection allows for an additional level of robustness to account for any poor/false structure detection during the initial mesh iterations where the mesh may be coarse or the approximated solution may inaccurately represent the active constraint arc.

\begin{table}[htp]
\centering
\caption{Estimated, $\hat{t}_s$, optimized, $t^*_s$, and relative difference of the activation and deactivation times obtained throughout mesh iterations involving structure detection for Example~$2$}
\begin{tabular}{| c | c c c | c c c | c |} 
 \hline
 \rule{0pt}{3ex}  & \multicolumn{3}{c|}{Activation time $[s]$} & \multicolumn{3}{c|}{Deactivation time $[s]$} & Stopping criteria\\[3pt] \hline
 \rule{0pt}{4ex} $M$ & $\hat{t}_{1}$ & $t^*_{1}$ &  $\dfrac{|\hat{t}-t^*|}{t^*} \times 100$ & $\hat{t}_{2}$ & $t^*_{2}$ &  $\dfrac{|\hat{t}-t^*|}{t^*} \times 100$ & $\max{\{\Delta c\}} \leq \epsilon_{c}$ \\[6pt] \hline \hline
 \rule{0pt}{3ex}1 &220.66 & 150.68 &  46.44 \% & 1012.96  & 1048.18 & 3.360 \%  & No \\[3pt] \hline 
 \rule{0pt}{3ex}2 &135.26 & 144.18 & 6.19 \%& 1062.67  & 1046.80 & 1.516 \% & No\\[3pt] \hline 
 \rule{0pt}{3ex}3 &138.85 & 139.52 & 0.48 \%& 1046.80  & 1046.75 & 0.005 \% & Yes\\[3pt] \hline 
\end{tabular}
\label{table:detection_performance}
\end{table}

\begin{figure}[htbp]
\centering
\vspace*{0.25cm}
\subfloat[Altitude, $h(t)$ vs. time, $t$\label{fig:altitude}]{\includegraphics[scale=0.4]{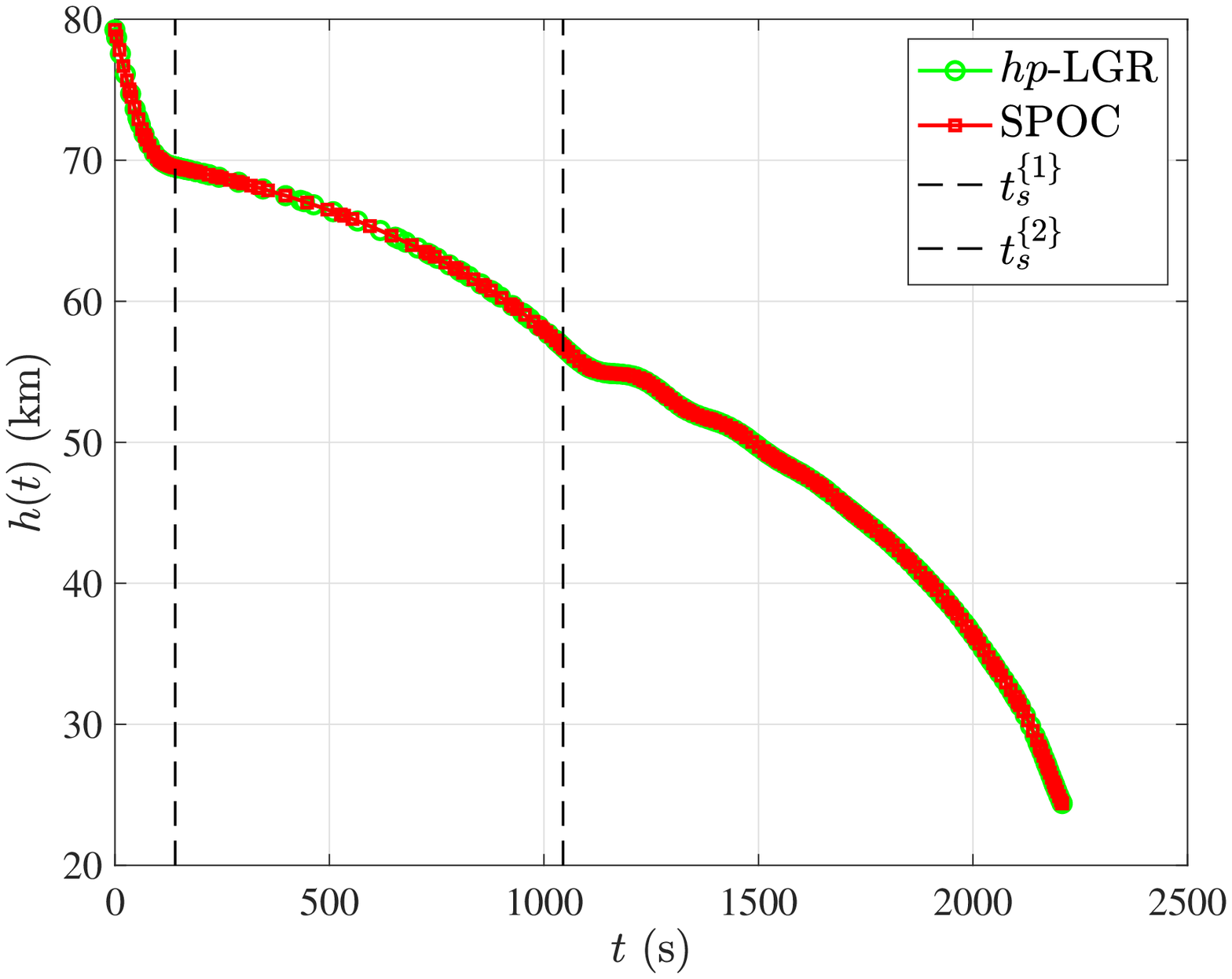}}~~\subfloat[Speed, $v(t)$ vs. time, $t$\label{fig:speed}]{\includegraphics[scale=0.4]{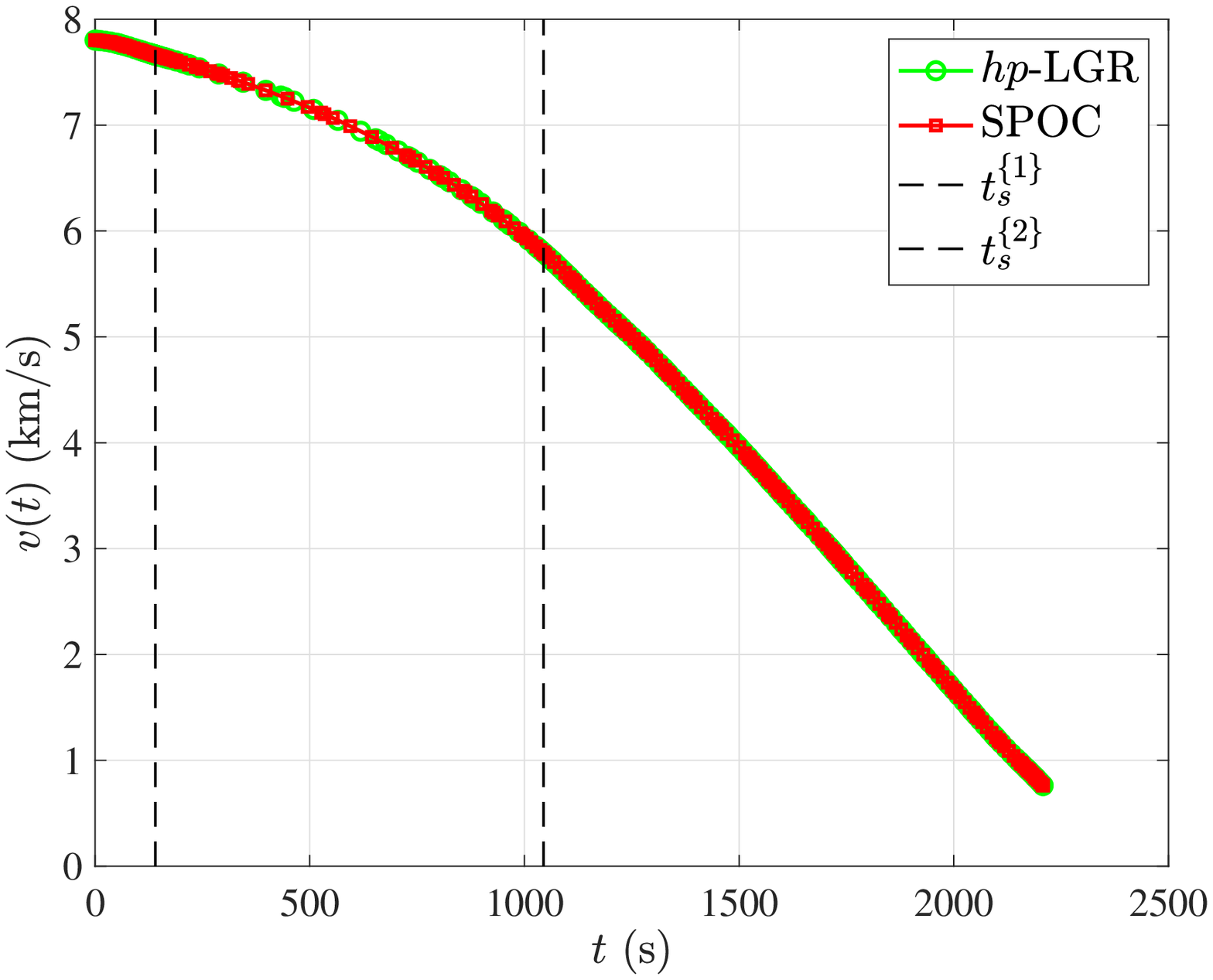}} \\ 
\subfloat[Latitude, $\phi(t)$ vs. longitude, $\theta(t)$\label{fig:latvlong}]{\includegraphics[scale=0.4]{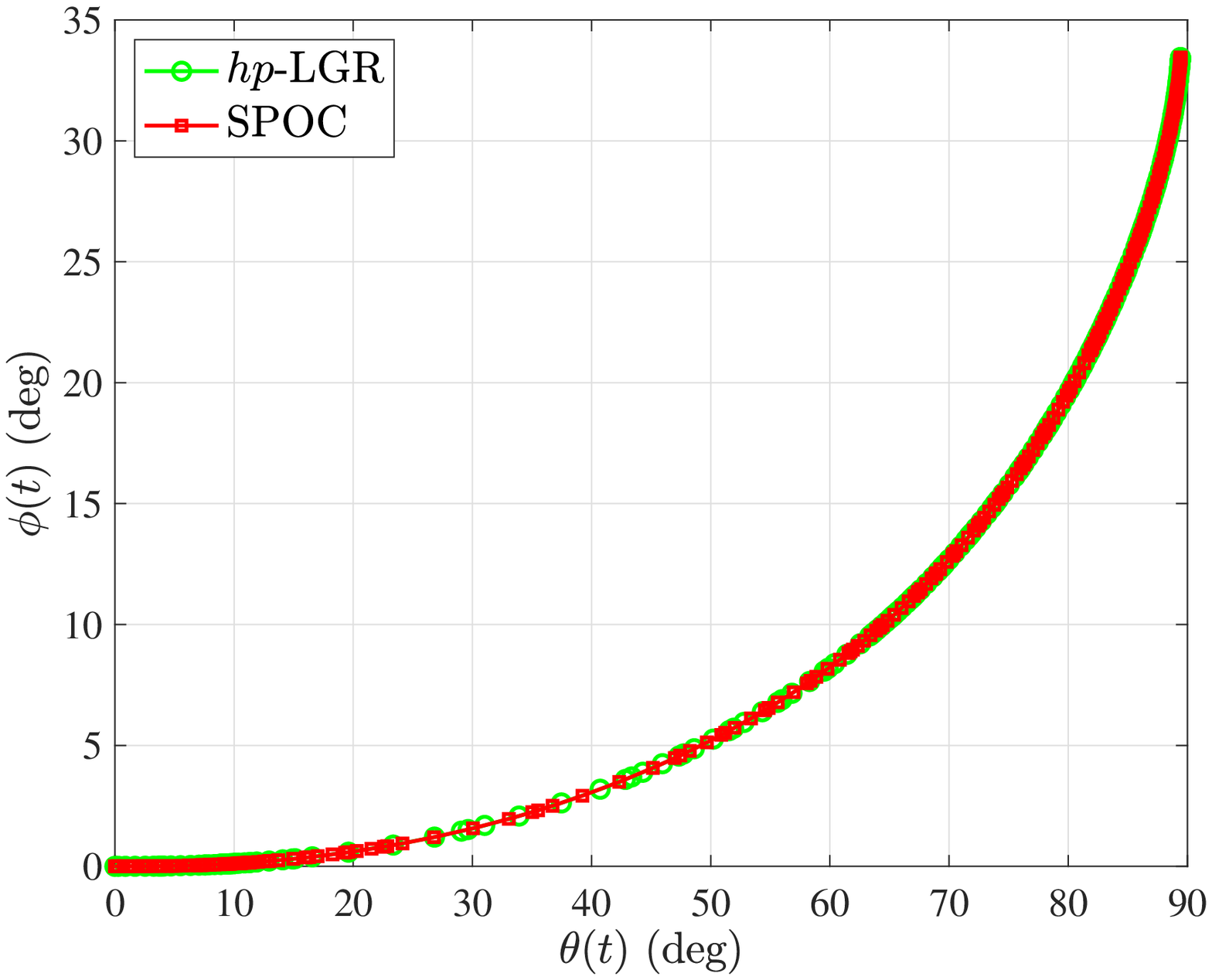}}~~\subfloat[Flight path angle, $\gamma(t)$, vs. time, $t$\label{fig:fpa}]{\includegraphics[scale=0.4]{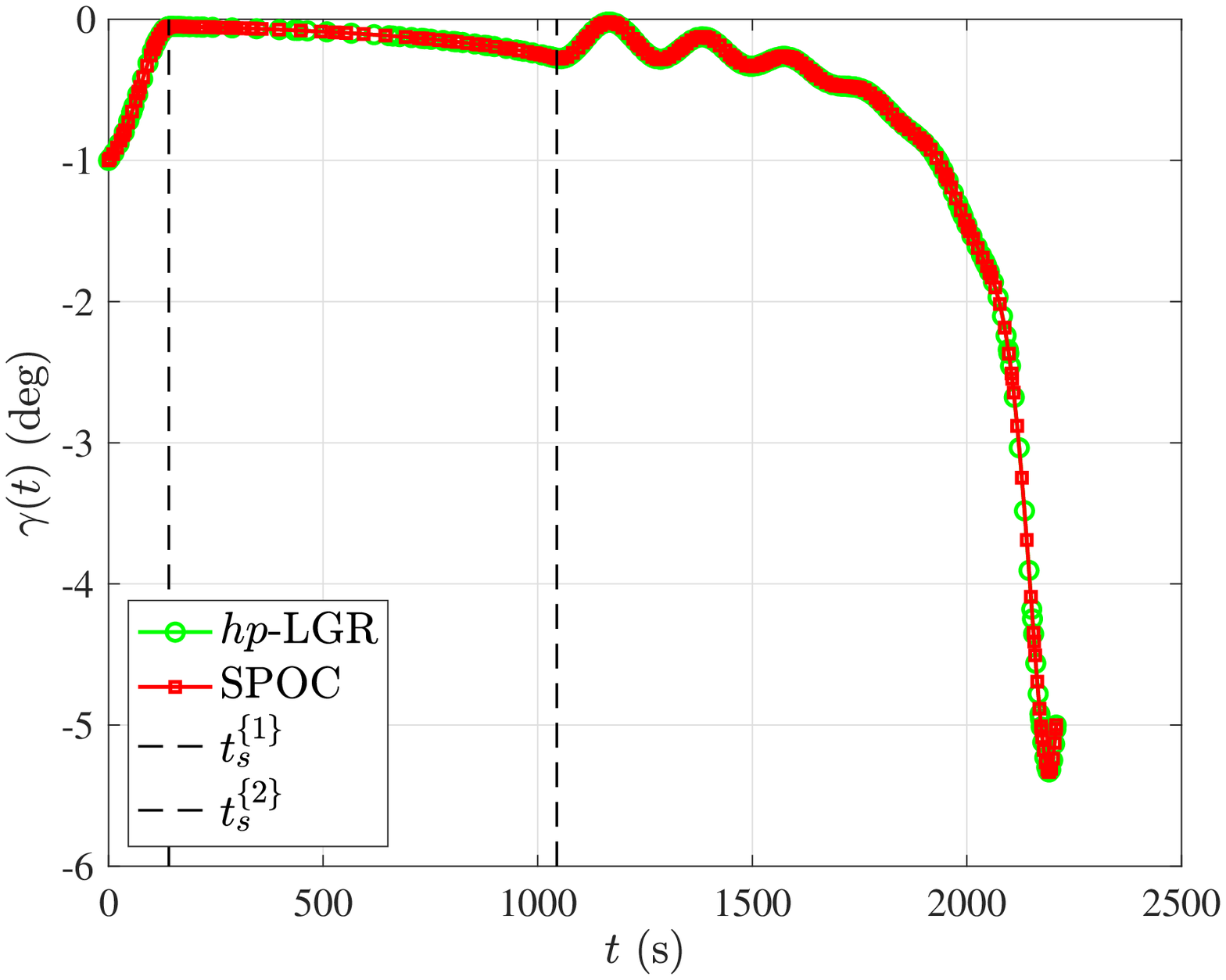}} \\ 
\subfloat[Angle of attack, $\alpha(t)$ vs. time, $t$\label{fig:aoa}]{\includegraphics[scale=0.4]{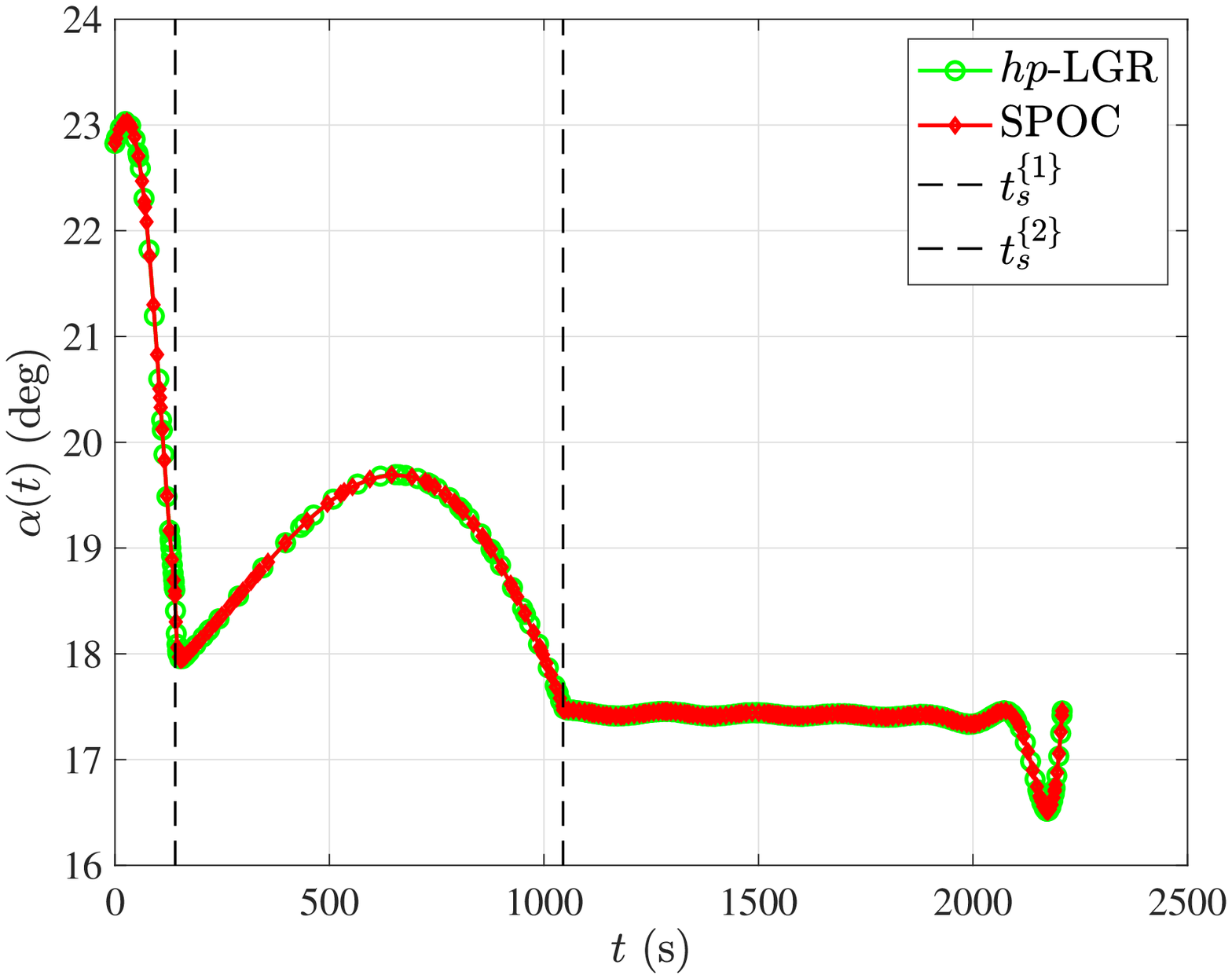}}~~\subfloat[Bank angle, $\sigma(t)$ vs. time, $t$\label{fig:bank_angle}]{\includegraphics[scale=0.4]{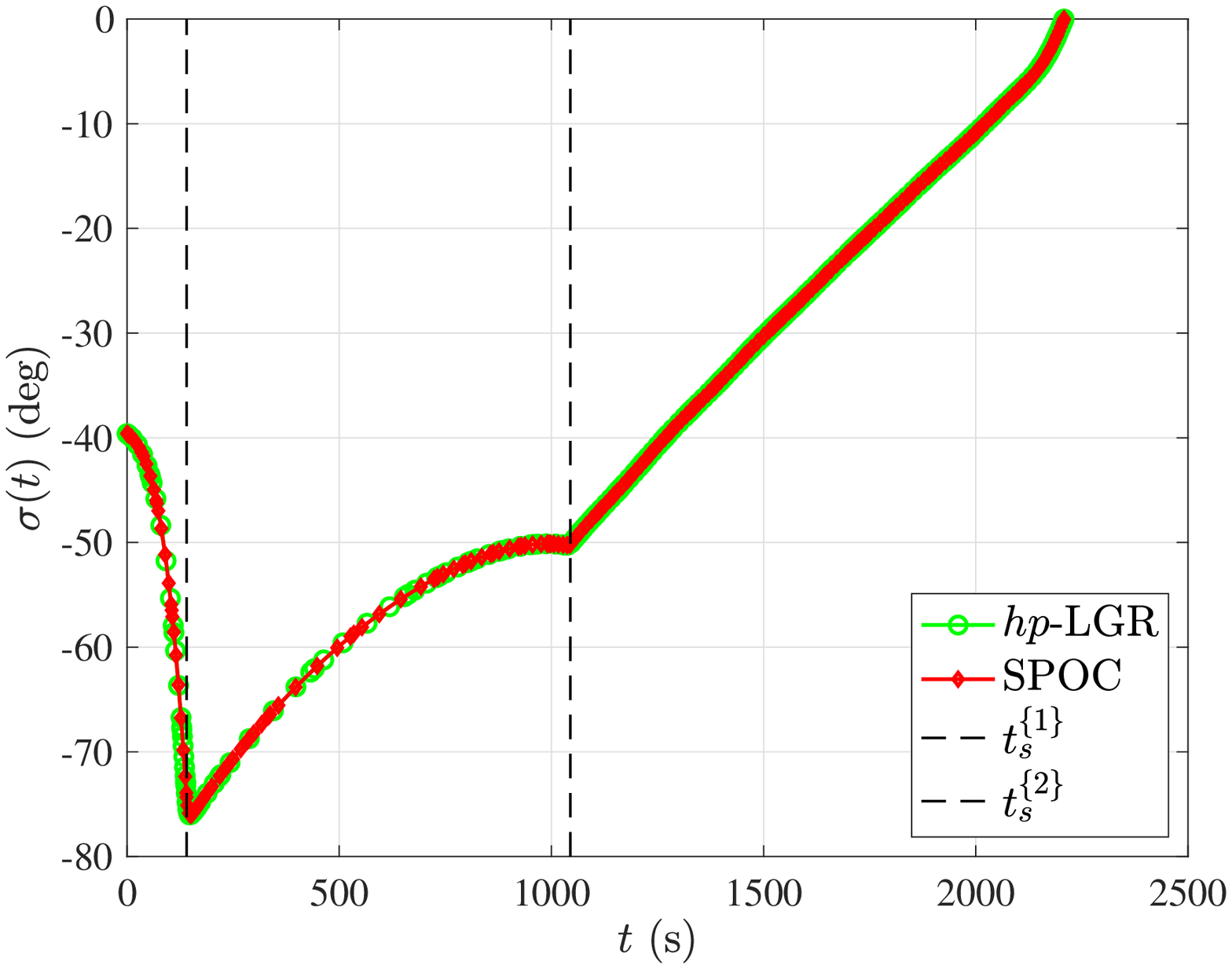}} \\ 
\caption{Solution to Example~$2$ for both the $hp$-LGR and SPOC methods, respectively.} 
\label{fig:RLVE_solution}
\end{figure}

For completeness, the state and control solutions obtained with both the $hp$-LGR and SPOC methods is provided in Figure~\ref{fig:RLVE_solution}. The solutions are found to be in close agreement with one another, which is expected given the computed objective values are in close agreement with one another. 

%%%%%%%%%%%%%%%%%%%%%%%%%%%%%%% SECTION 7 %%%%%%%%%%%%%%%%%%%%%%%%%%%%%%%%%%%%%%%%%%%%%%%%%%%%%%%
\section{Limitations of the SPOC Method}\label{section:limitations}
Like most computational methods, some trade-offs and limitations exist with the SPOC method. The most apparent trade-off that exists within the SPOC method is between accuracy and computational time. Both examples presented in this work display an improvement in accuracy as compared against the $hp$-LGR method, but at the expense of a larger computation time. One such limitation to the SPOC method is the assumption that the functions defining the state-variable inequality path constraints are continuous and differentiable, which may not always be the case. For example, one might have a set of tabular data that defines the lower/upper limits on a particular state-variable inequality path constraint. If this is the case, the SPOC method would require the user to perform some form of a function approximation leading to a potential loss of accuracy in capturing the true behavior of the constraint. \par
Another limitation to the SPOC method is the dependence on an initial solution and a user specified detection tolerance $\epsilon$. Specifically, the first pass through the structure detection algorithm is dependent on the quality of the initial solution and how "tight" the detection tolerance is set. For example, decreasing the detection tolerance tightens the window for what defines an activation and deactivation point, but may lead to false activation and deactivation points if the quality of the provided solution on the previous mesh contains chattering. Conversely, increasing the detection tolerance would loosen the window for what defines an activation and deactivation point, but may lead to detecting activation and deactivation points far from the true points. While the default detection tolerance given in Sect.~\ref{section:structure_detect} worked for the two examples shown, it may require tuning from a user for more complex problems. \par
Finally, while results obtained with the SPOC method were compared against the $hp$-LGR method, it also uses the $hp$-LGR method within to suggest new meshes when the maximum relative error hasn't been met. A potential direction for future work could investigate a refinement algorithm tailored towards problems with state-variable inequality path constraints, and even problems with other types of constraints such as control constraints and mixed-state control constraints. 

%%%%%%%%%%%%%%%%%%%%%%%%%%%%%%% SECTION 8 %%%%%%%%%%%%%%%%%%%%%%%%%%%%%%%%%%%%%%%%%%%%%%%%%%%%%%%
\section{Conclusions}\label{section:conclusions}
A method for solving optimal control problems with state-variable inequality path constraints has been developed. The method implements a detection algorithm that returns estimated activation and deactivation times of the active state-variable inequality path constraints. Based on the detection results, the method partitions the original time domain into multiple subdomains where the constraints are deemed to either be active or inactive. Within the multiple-domain formulation, the method introduces additional decision variables in the optimization process which represent the estimated activation and deactivation times, thus allowing for the exact locations of the activation and deactivation times to be optimized. For domains categorized as constrained, the lowest time derivative of the active path constraints which are explicit in a control variable is set to zero, while the preceding derivatives and the constraints themselves are set to zero at the beginning of the constrained domain. Lastly, two stopping criteria are checked to determine if an adaptive mesh refinement procedure must be employed and/or the structure detection algorithm must be used to re-detect the active constraints. The method has been demonstrated on two examples from the existing literature and compared against the method of Ref. \cite{Liu2018}. It is found that the method is capable of producing an accurate solution to optimal control problems involving state-variable inequality path constraints.
\section{Acknowledgments}
The authors gratefully acknowledge support for this research from the U.S. National Science Foundation under grant CMMI-2031213, the U.S. Office of Naval Research under grant N00014-22-1-2397, and from the U.S. Air Force Research Laboratory under contract FA8651-21-F-1041.
\renewcommand{\baselinestretch}{1}\normalsize\normalfont
%\bibliography{References}
%\bibliographystyle{ieeetr}

\end{document}